\documentclass[11pt,oneside]{amsart}
\usepackage{amsmath}
\usepackage{amssymb}
\usepackage{amscd}
\usepackage{amsthm}
\usepackage{latexsym}
\usepackage[pdftex]{graphicx}
\usepackage{epstopdf}
\usepackage{pinlabel}
\usepackage{enumerate}
\usepackage{epstopdf}

\numberwithin{equation}{section}
\theoremstyle{plain}

\newtheorem{thm}{Theorem}[section]

\newtheorem{thm-defn}[thm]{Theorem/Definition}

\newtheorem{cor}[thm]{Corollary}

\theoremstyle{definition}
\newtheorem{defn}[thm]{Definition}

\theoremstyle{remark}

\newcommand{\C}{\mathbb C}

\newcommand{\Q}{\mathbb Q}

\newcommand{\Z}{\mathbb Z}

\begin{document}

\begin{abstract}
The rational ho\-mol\-ogy balls $B_n$ ap\-peared in Fin\-tu\-shel and Stern's ration\-al blow-down construction \cite{FSrbd} and were subsequently used in \cite{Pa2,FS4} to construct exotic smooth manifolds with small Euler numbers. We show that a large class of smooth $4$-manifolds have all of the $B_n$'s for odd $n \geq 3$ embedded in them. In particular, we give explicit examples, using Kirby calculus, of several families of smooth embeddings of the rational homology balls $B_n$.
\end{abstract}

\title[Smooth Embeddings of Rational Homology Balls]{Smooth Embeddings of Rational Homology Balls}
\author{Tatyana Khodorovskiy}

\maketitle

\section{Introduction}
\begin{figure}[ht!]
\labellist
\small\hair 2pt
\pinlabel $-(n+2)$ at 30 4.7
\pinlabel $-2$ at 60 4.7
\pinlabel $-2$ at 90 4.7
\pinlabel $-2$ at 173 4.7
\pinlabel $-2$ at 203 4.7
\pinlabel $S_1$ at 32 2.5 
\pinlabel $S_2$ at 62 2.5 
\pinlabel $S_3$ at 92 2.5 
\pinlabel $S_{n-2}$ at 175 2.5 
\pinlabel $S_{n-1}$ at 205 2.5 
\endlabellist
\centering
\includegraphics[width=120mm,height=30mm]{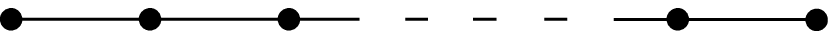}
\caption{{\bf Plumbing diagram of $C_n$, $n\geq2$}}
\label{f:cn}\end{figure}

In 1997, Fintushel and Stern \cite{FSrbd} defined the rational blow-down operation for smooth $4$-manifolds, a generalization of the standard blow-down operation. For smooth $4$-manifolds, the standard blow-down is performed by removing a neighborhood of a sphere with self-intersection $(-1)$ and replacing it with a standard $4$-ball $B^4$. The rational blow-down involves replacing a negative definite plumbing $4$-manifold with a rational homology ball. In order to define it, we first begin with a description of the negative definite plumbing $4$-manifold $C_n$, $n \geq 2$, as seen in Figure~\ref{f:cn}, where each dot represents a sphere, $S_i$, in the plumbing configuration. The integers above the dots are the self-intersection numbers of the plumbed spheres: $[S_1]^2 = -(n+2)$ and $[S_i]^2 = -2$ for $2 \leq i \leq n-1$.

The boundary of $C_n$ is the lens space $L(n^2,n-1)$, thus $\pi_1(\partial C_n) \cong H_1(\partial C_n ; \Z) \cong \Z/n^2\Z$. (Note, when we write the lens space $L(p,q)$, we mean it is the $3$-manifold obtained by performing $-\frac{p}{q}$ surgery on the unknot.) This follows from the fact that $[-n-2, -2, \ldots -2]$, with $(n-2)$ many $(-2)$'s is the continued fraction expansion of $\frac{n^2}{1-n}$.

\begin{figure}[ht!]
\labellist
\small\hair 2pt
\pinlabel $n-1$ at -25 240
\pinlabel $n$ at 235 160
\endlabellist
\centering
\includegraphics[scale=0.25]{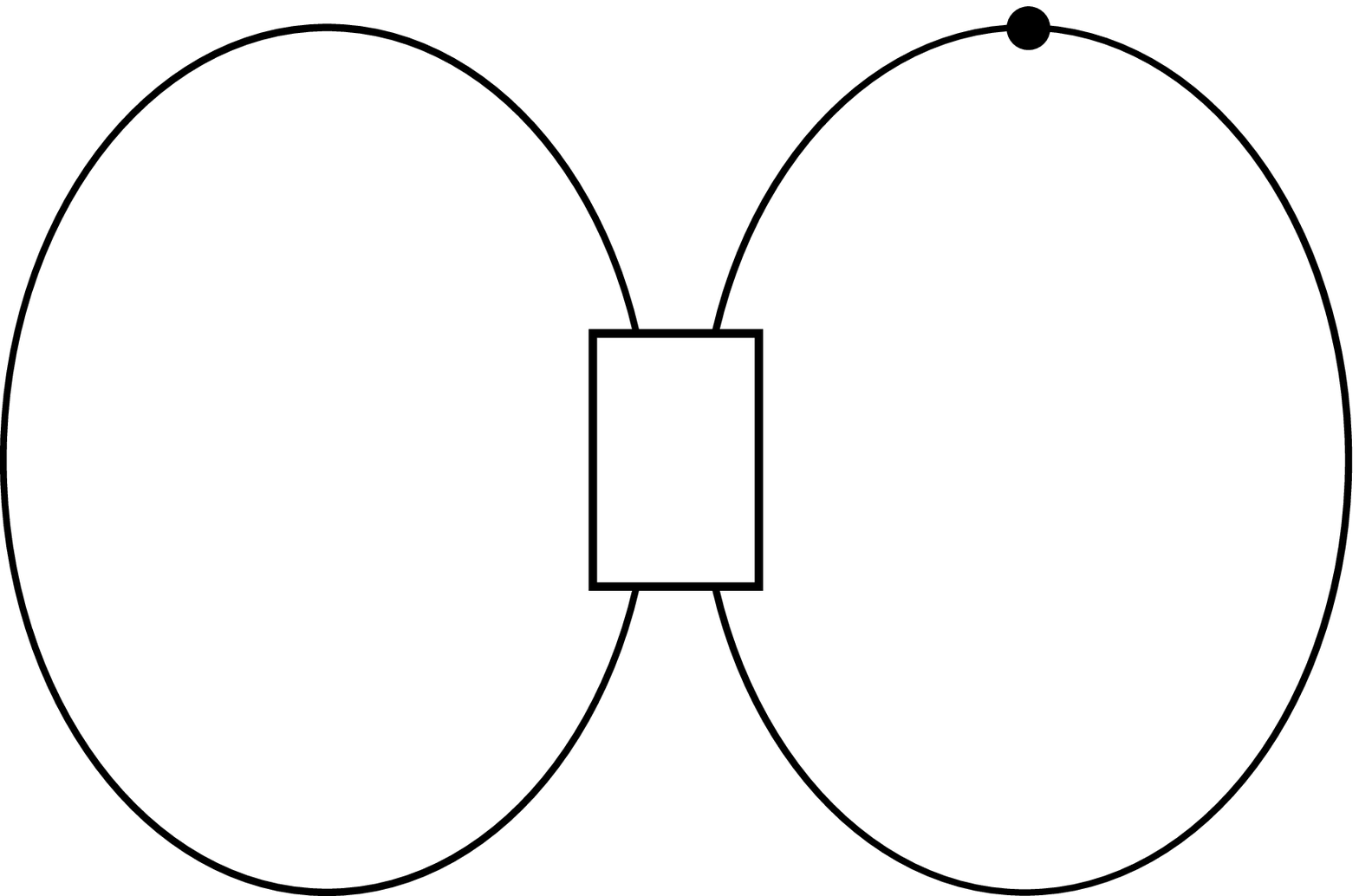}
\caption{{\bf Kirby diagram of $B_n$}}
\label{f:bn}
\end{figure}

Let $B_n$ be the $4$-manifold as defined by the Kirby diagram in Figure~\ref{f:bn} (for a more extensive description of $B_n$, see section~\ref{sec:bn}). The manifold $B_n$ is a rational homology ball, i.e. $H_*(B_n;\Q) \cong H_*(B^4;\Q)$. The boundary of $B_n$ is also the lens space $L(n^2,n-1)$ \cite{CaHa}. Moreover, any self-diffeomorphism of $\partial B_n$ extends to $B_n$ \cite{FSrbd}. Now, we can define the rational blow-down of a $4$-manifold $X$ below in Definition~\ref{def:rbd}. Fintushel and Stern \cite{FSrbd} also showed how to compute Seiberg-Witten and Donaldson invariants of $X_{(n)}$ from the respective invariants of $X$. In addition, they showed that certain smooth logarithmic transforms can be alternatively expressed as a series of blow-ups and rational blow-downs. The rational blow-down was particularly useful in constructing $4$-manifolds with exotic smooth structures, with small Euler numbers.

\begin{defn}
\label{def:rbd}
(\cite{FSrbd}, also see \cite{GS}) Let $X$ be a smooth $4$-manifold. Assume that $C_n$ embeds in $X$, so that $X = C_n \cup_{L(n^2,n-1)} X_0$. The 4-manifold $X_{(n)} = B_n \cup_{L(n^2,n-1)} X_0$ is by definition the \textit{rational blow-down} of $X$ along the given copy of $C_n$.
\end{defn}

One can define the (smooth) \textit{rational blow-up} operation in a similar manner: if there exists a smoothly embedded $B_n$ in a $4$-manifold $X$ then one can rationally blow-up $X$ by removing the $B_n$ and gluing in the $C_n$, along the lens space $L(n^2,n-1)$. Consequently, one can ask: which $4$-manifolds can be smoothly rationally blown up? Equivalently, which $4$-manifolds contain a smoothly embedded rational homology ball $B_n$? We prove the following results regarding smooth embeddings of the rational homology balls $B_n$:

\begin{thm}
\label{thm:smoothbn1}
Let $V_{-n-1}$ be a neighborhood of a sphere with self-intersection number $(-n-1)$. There exists an embedding of the rational homology balls $B_n \hookrightarrow V_{-n-1}$, for all $n \geq 2$.
\end{thm}

\begin{thm}
\label{thm:smoothbn2}
Let $V_{-4}$ be a neighborhood of a sphere with self-intersection number $(-4)$. For all $n \geq 3$ odd, there exists an embedding of the rational homology balls $B_n \hookrightarrow V_{-4}$. For all $n \geq 2$ even, there exists an embedding of the rational homology balls $B_n \hookrightarrow B_2 \# \overline{\C P^2}$.
\end{thm}

\noindent Theorems~\ref{thm:smoothbn1} and~\ref{thm:smoothbn2} above show that there is little obstruction to \textit{smoothly} embedding the rational homology balls $B_n$ into a smooth $4$-manifold. In particular, Theorem~\ref{thm:smoothbn2} implies that if a smooth $4$-manifold $X$ contains a sphere with self-intersection $(-4)$, then one can smoothly embed the rational homology balls $B_n$ into $X$ for all odd $n \geq 3$. 

One of the implications of Theorem~\ref{thm:smoothbn2} is that for a given smooth $4$-manifold $X$, there does not exist an $N$, such that for all $n \geq N$ one cannot find a smooth embedding $B_n \hookrightarrow X$. In the setting of this sort in algebraic geometry, for rational homology ball smoothings of certain surface singularities, such a bound on $n$ does exist, in terms of $(c_1^2,\chi_h)$ invariants of an algebraic surface \cite{KoSB,Wa}.

In section~\ref{sec:bn} we give a brief description of the rational homology balls $B_n$. In section~\ref{sec:smoothemb} we describe some previously known embeddings of $B_n$ in order to illustrate the differences between them and those embeddings in Theorems~\ref{thm:smoothbn1} and~\ref{thm:smoothbn2}. In sections~\ref{sec:pfthm1} and~\ref{sec:pfthm2} we prove Theorems~\ref{thm:smoothbn1} and~\ref{thm:smoothbn2} respectively. Finally, in section~\ref{sec:simple} we define the notion of ``simple" embeddings of $B_n$.

\section{Description of the rational homology balls $B_n$}
\label{sec:bn}

\begin{figure}[ht!]
\labellist
\small\hair 2pt
\pinlabel $n-1$ at 475 18
\pinlabel $\}n$ at 190 65
\endlabellist
\centering
\includegraphics[scale=0.5]{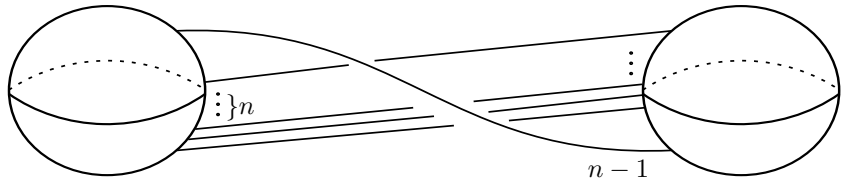}
\caption{{\bf Another Kirby diagram of $B_n$}}
\label{f:bnsph}
\end{figure}

There are several ways to give a description of the rational homology balls $B_n$. One of them is a Kirby calculus diagram seen in Figure~\ref{f:bn}. This represents the following handle decomposition: Start with a 0-handle, a standard 4-disk $D^4$, attach to it a 1-handle $D^1 \times D^3$. Call the resultant space $X_1$, it is diffeomorphic to $S^1 \times D^3$ and has boundary $\partial X_1 = S^1 \times S^2$. Finally, we attach a 2-handle $D^2 \times D^2$. The boundary of the core disk of the 2-handle gets attached to the closed curve, $K$, in $\partial X_1$ which wraps $n$ times around the $S^1 \times \ast$ in $S^1 \times S^2$. We can also represent $B_n$ by a slightly different Kirby diagram, which is more cumbersome to manipulate but is more visually informative, as seen in Figure~\ref{f:bnsph}, where the 1-handle is represented by a pair of balls. It is worthwhile to note that $B_2$ can also be described as an unoriented disk bundle over $\mathbb{RP}^2$, where $K$ is the boundary of the Mobius band in $\mathbb{RP}^2$.

\section{Embeddings of $B_n$ obtained from blow-ups of elliptical surfaces}
\label{sec:smoothemb}

\begin{figure}[ht!]
\labellist
\small\hair 2pt
\pinlabel $0$ at 27 45
\pinlabel $-4$ at 70 40
\pinlabel $-1$ at 110 35
\pinlabel $-5$ at 165 35
\pinlabel $-2$ at 191 35
\pinlabel $-1$ at 175 20
\pinlabel $-n-2$ at 310 48
\pinlabel $-1$ at 303 35
\pinlabel $-2$ at 310 20
\pinlabel $-2$ at 346 47
\pinlabel $-2$ at 364 35
\endlabellist
\centering
\includegraphics[height=40mm, width=125mm]{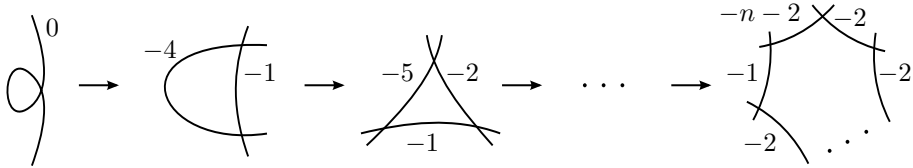}
\caption{Fishtail fiber in $E(m)$ blown up $(n-1)$ times}
\label{f:fishfiber}
\end{figure}

In the existing literature, the straightforward examples of smooth $4$-manifolds containing smoothly embedded rational homology balls $B_n$, are those manifolds obtained via a rational blow-down. Examples of such manifolds first appeared in Fintushel and Stern's original paper \cite{FSrbd} on rational blow-downs: logarithmic transforms $E(m)_n$ of elliptic surfaces $E(m)$. In these manifolds, one starts with a fishtail fiber of $E(m)$, which has homological self-intersection $0$, blows it up $(n-2)$ times, and then one obtains a configuration of spheres $C_n$, which one rationally blows down (see Figure~\ref{f:fishfiber}). Consequently, one obtains a manifold $E(m)_n$, having the same $(c_1^2,c_2)$ numbers but different Seiberg-Witten invariants as $E(m)$, which contains an embedded rational homology ball $B_n$.

Most other examples of smooth $4$-manifolds that contain embedded rational homology balls $B_n$, are obtained in a similar manner, one blows up a smooth manifold several times, then finds a particular configuration of spheres $C_n$ which one rationally blows down. Often, one ends up with a manifold with lower betti number $b_2^-$ than the original manifold one started with. In fact, in a lot of these examples, since one can compute the betti numbers of the resultant manifold, by Freedman's theorem \cite{Fr,FrQu} one can conclude they are homeomorphic to $k\C P^2 \# \ell \overline{\C P^2}$, for some $k$ and $\ell$. However, after a computation of the Seiberg-Witten invariants, one can often show that the resultant manifolds are not diffeomorphic to $k\C P^2 \# \ell \overline{\C P^2}$, and thus possess an exotic smooth structure, which is frequently the goal. In fact, one can sometimes find an infinite family of exotic $4$-manifolds which are homeomorphic but not diffeomorphic to $k\C P^2 \# \ell \overline{\C P^2}$. For example, using these techniques, exotic $\C P^2 \# 7\overline{\C P^2}$ manifolds were constructed in \cite{Pa2}. Additionally, using a generalized rational blow-down \cite{Park}, exotic  $\C P^2 \# 6\overline{\C P^2}$ manifolds were constructed in \cite{StSz}.

\section{Proof of Theorem~\ref{thm:smoothbn1}}
\label{sec:pfthm1}

In this section we prove Theorem~\ref{thm:smoothbn2}.

\begin{proof}{\textit{Proof of Theorem~\ref{thm:smoothbn1}}.}

\begin{figure}[ht]
\begin{minipage}[b]{0.55\linewidth}
\centering
\includegraphics[scale=0.22]{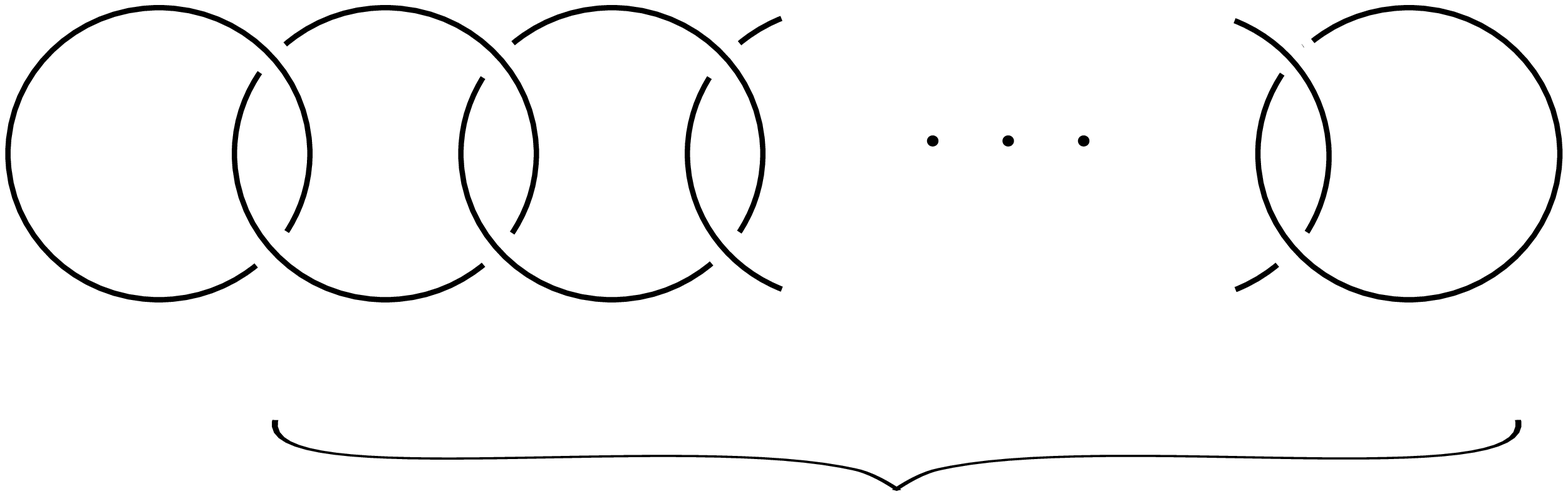}
\caption{ }
\labellist
\small\hair 2pt
\pinlabel $-n-2$ at -720 210
\pinlabel $-2$ at -580 210
\pinlabel $-2$ at -490 210
\pinlabel $-2$ at -190 210
\pinlabel $n-2$ at -410 -10
\pinlabel $-1$ at 185 170
\pinlabel $-1$ at 175 140
\pinlabel $-1$ at 180 65
\pinlabel $0$ at 430 230
\pinlabel $n-1$ at 18 120
\pinlabel $n$ at 288 115
\pinlabel $S_1$ at -690 60
\pinlabel $S_2$ at -600 60
\pinlabel $S_3$ at -510 60
\pinlabel $S_{n-1}$ at -190 60
\endlabellist
\label{f:KirbyCnCircles1}
\end{minipage}
\hspace{0.5cm}
\begin{minipage}[b]{0.35\linewidth}
\centering
\includegraphics[scale=0.28]{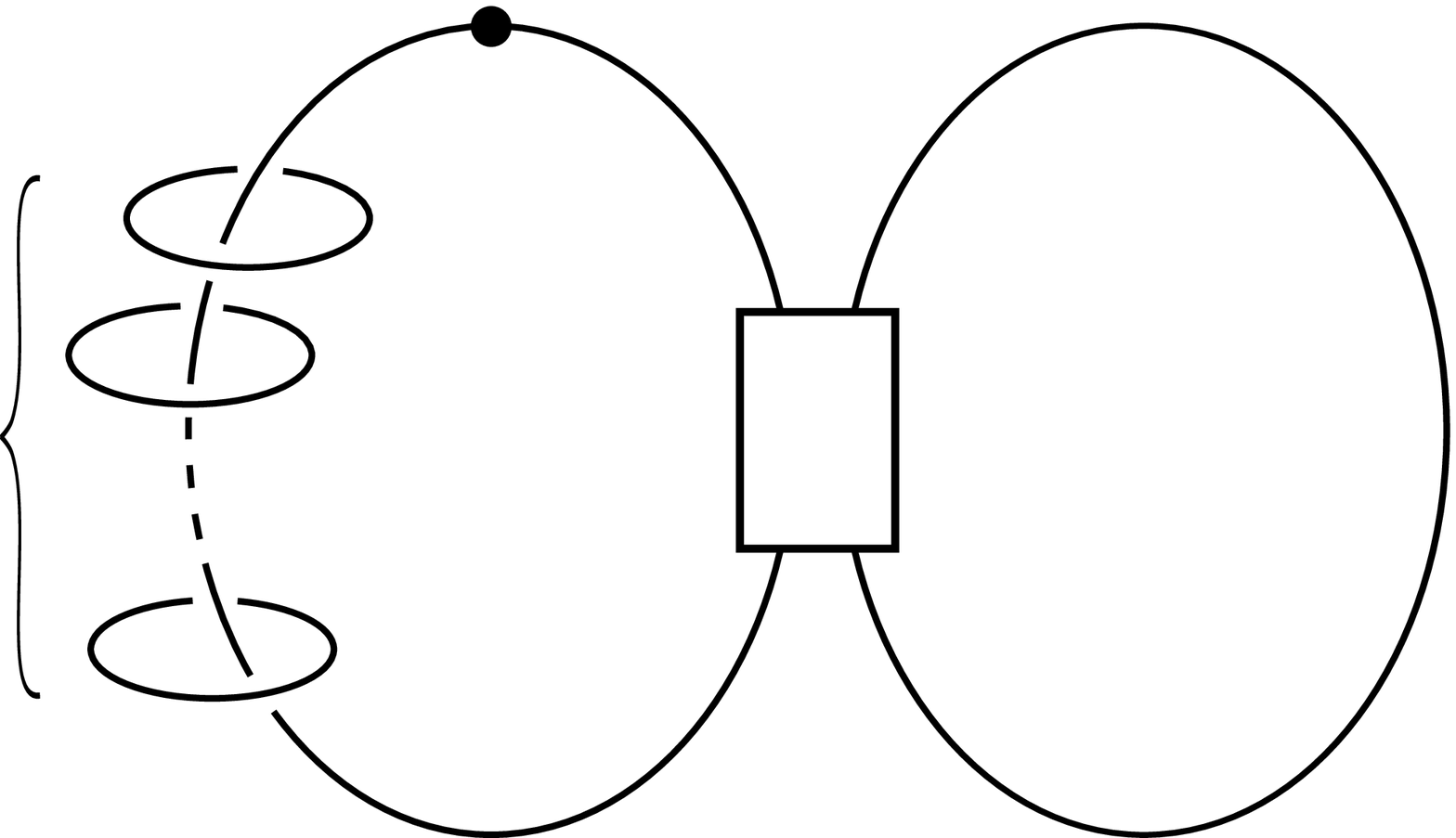}
\caption{ }
\label{f:KirbyCnCircles2}
\end{minipage}
\end{figure}

We prove this theorem using Kirby calculus (see \cite{GS} for detailed exposition). We start with the Kirby diagram for $V_{-n-1}$, Figure~\ref{f:KirbyThmA1}, blow it up $(n-1)$ times, and obtain the configuration of spheres $C_n$ with an additional sphere $\Sigma_{-1}$ with  $[\Sigma_{-1}]^2 = -1$, attached to the last sphere with self-intersection $(-2)$, $S_{n-1}$, Figure~\ref{f:KirbyThmA4}. In Figures~\ref{f:KirbyThmA5}-\ref{f:KirbyThmA9} we proceed to do the standard Kirby calculus manipulation where one changes the Kirby diagram of $C_n$ from the one in Figure~\ref{f:KirbyCnCircles1} to the one in Figure~\ref{f:KirbyCnCircles2}, (see \cite{GS}, p. 516), by first adding a cancelling $1/2$ handle pair (Figure~\ref{f:KirbyThmA5}) and performing a series of handleslides. However, in our case the additional sphere $\Sigma_{-1}$ is present, and intersects with the $C_n$ configuration in a non-trivial way. As a result, when we perform the last handleslide to get $C_n$ to look like Figure~\ref{f:KirbyCnCircles2}, the sphere $\Sigma_{-1}$ intersects with $C_n$ as seen in Figure~\ref{f:KirbyThmA9}. 

\begin{figure}[ht]
\begin{minipage}[b]{0.30\linewidth}
\centering
\includegraphics[scale=0.28]{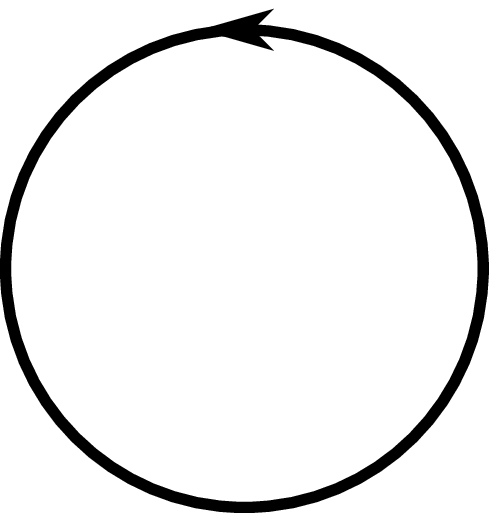}
\caption{$V_{-n-1}$}
\labellist
\small\hair 2pt
\pinlabel $-n-1$ at -260 135
\pinlabel $-n-2$ at 40 135
\pinlabel $-1$ at 200 130
\pinlabel $-n-2$ at 350 135
\pinlabel $-2$ at 500 130
\pinlabel $-1$ at 590 130
\endlabellist
\label{f:KirbyThmA1}
\end{minipage}
\begin{minipage}[b]{0.33\linewidth}
\centering
\includegraphics[height=20mm, width=30mm]{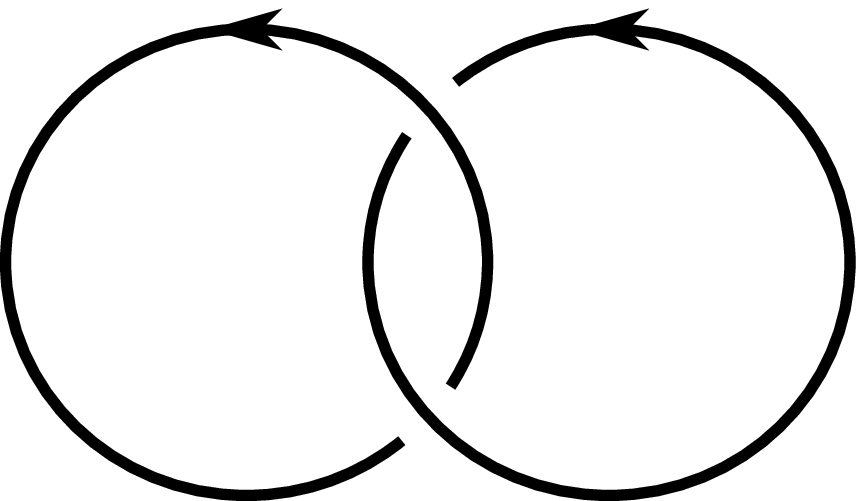}
\caption{$V_{-n-1} \# \overline{\C P^2}$ }
\label{f:KirbyThmA2}
\end{minipage}
\begin{minipage}[b]{0.34\linewidth}
\centering
\includegraphics[scale=0.28]{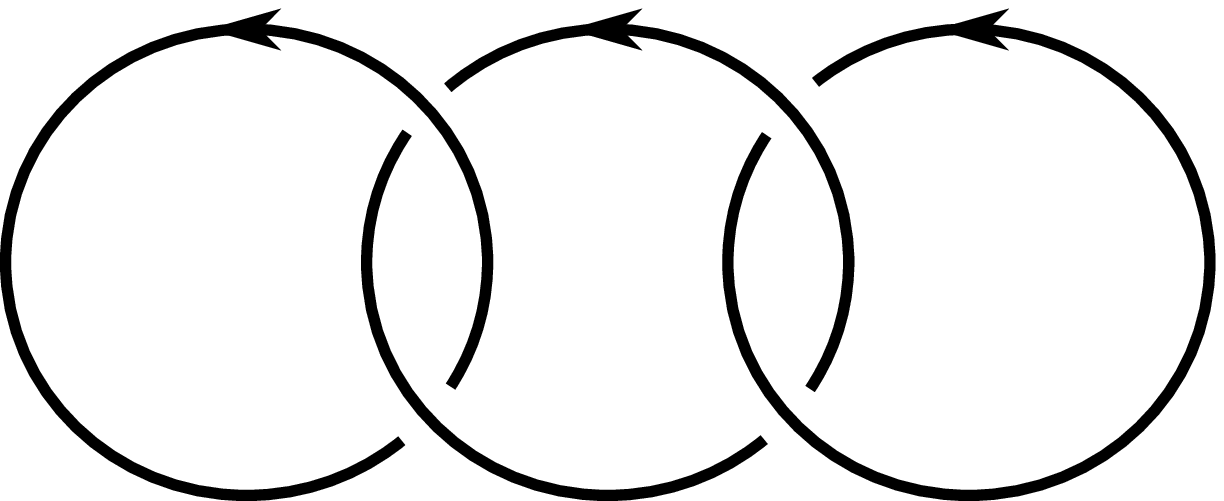}
\caption{$V_{-n-1} \# 2\overline{\C P^2}$ }
\label{f:KirbyThmA3}
\end{minipage}
\end{figure}

\begin{figure}[ht!]
\labellist
\small\hair 2pt
\pinlabel $-n-2$ at 40 175
\pinlabel $-2$ at 230 170
\pinlabel $-2$ at 330 170
\pinlabel $-2$ at 660 170
\pinlabel $-1$ at 770 170
\pinlabel $n-2$ at 405 -5
\endlabellist
\centering
\includegraphics[height=25mm, width=90mm]{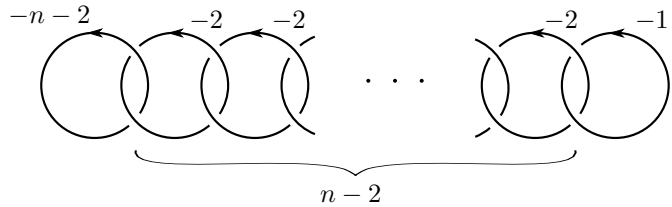}
\caption{$V_{-n-1} \# (n-1)\overline{\C P^2}$ }
\label{f:KirbyThmA4}
\end{figure}

\begin{figure}[ht!]
\labellist
\small\hair 2pt
\pinlabel $-n-2$ at 220 35
\pinlabel $-2$ at 370 70
\pinlabel $-2$ at 390 95
\pinlabel $-2$ at 550 95
\pinlabel $-1$ at 600 95
\pinlabel $n-2$ at 430 175
\pinlabel $-1$ at 250 160
\endlabellist
\centering
\includegraphics[height=45mm, width=125mm]{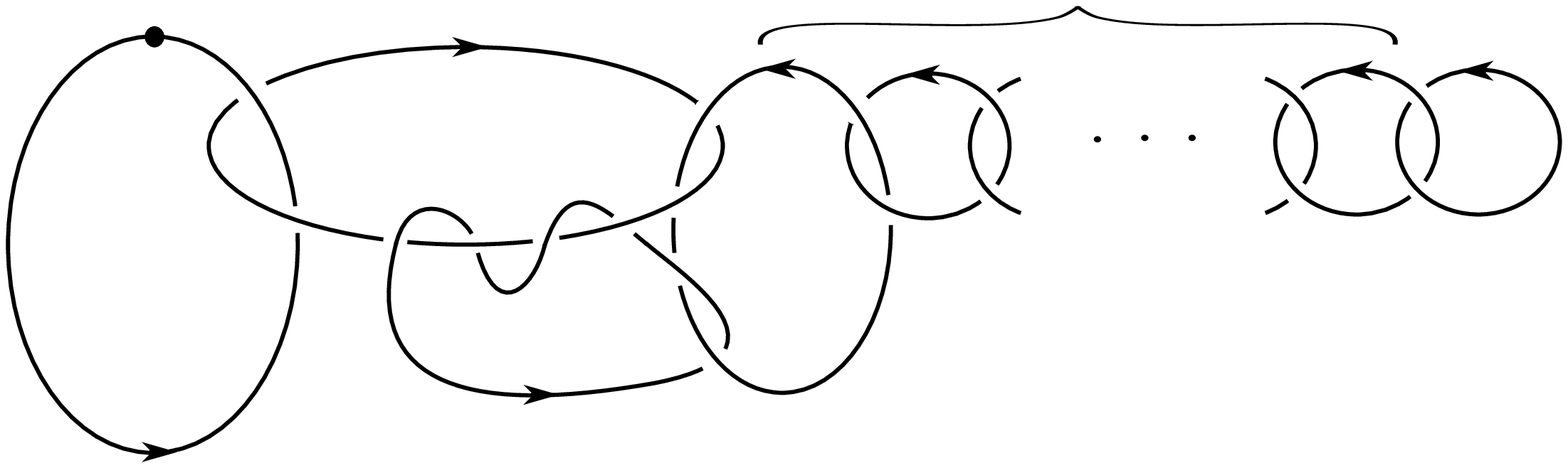}
\caption{ }
\label{f:KirbyThmA5}
\end{figure}

\begin{figure}[ht!]
\labellist
\small\hair 2pt
\pinlabel $-n+1$ at 570 60
\pinlabel $-2$ at 275 87
\pinlabel $-2$ at 325 87
\pinlabel $-2$ at 470 87
\pinlabel $-1$ at 520 87
\pinlabel $n-2$ at 370 55
\pinlabel $-1$ at 220 153
\endlabellist
\centering
\includegraphics[height=45mm, width=120mm]{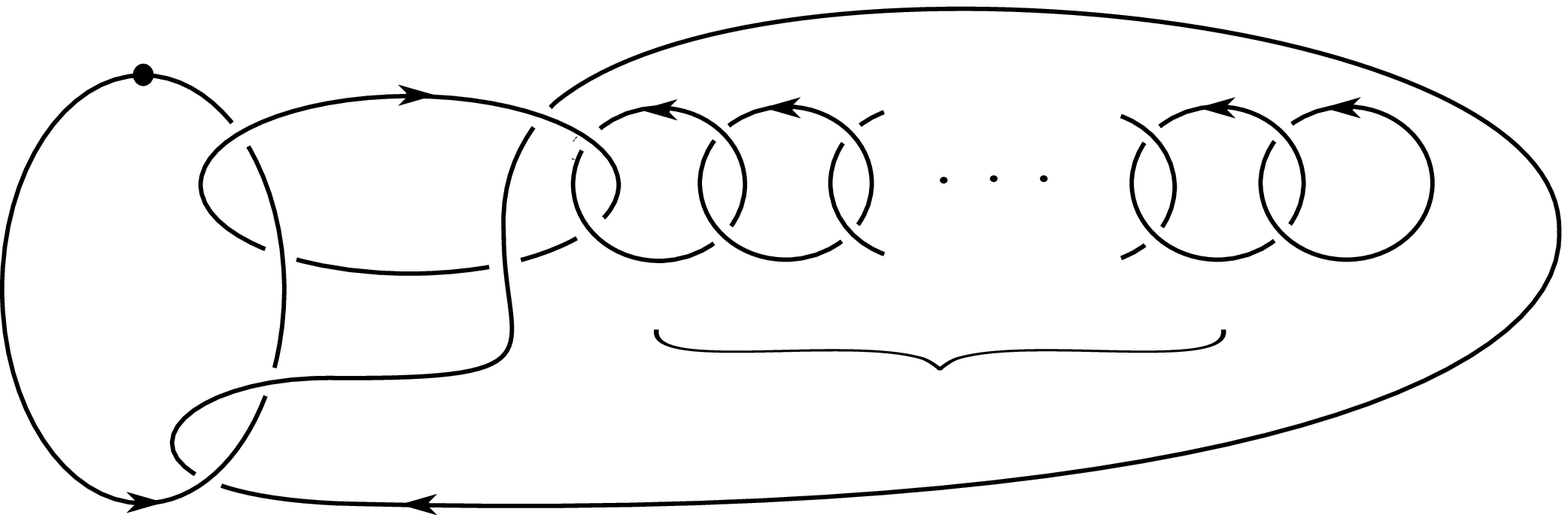}
\caption{ }
\label{f:KirbyThmA6}
\end{figure}

\begin{figure}[ht]
\begin{minipage}[b]{0.45\linewidth}
\centering
\includegraphics[scale=0.47]{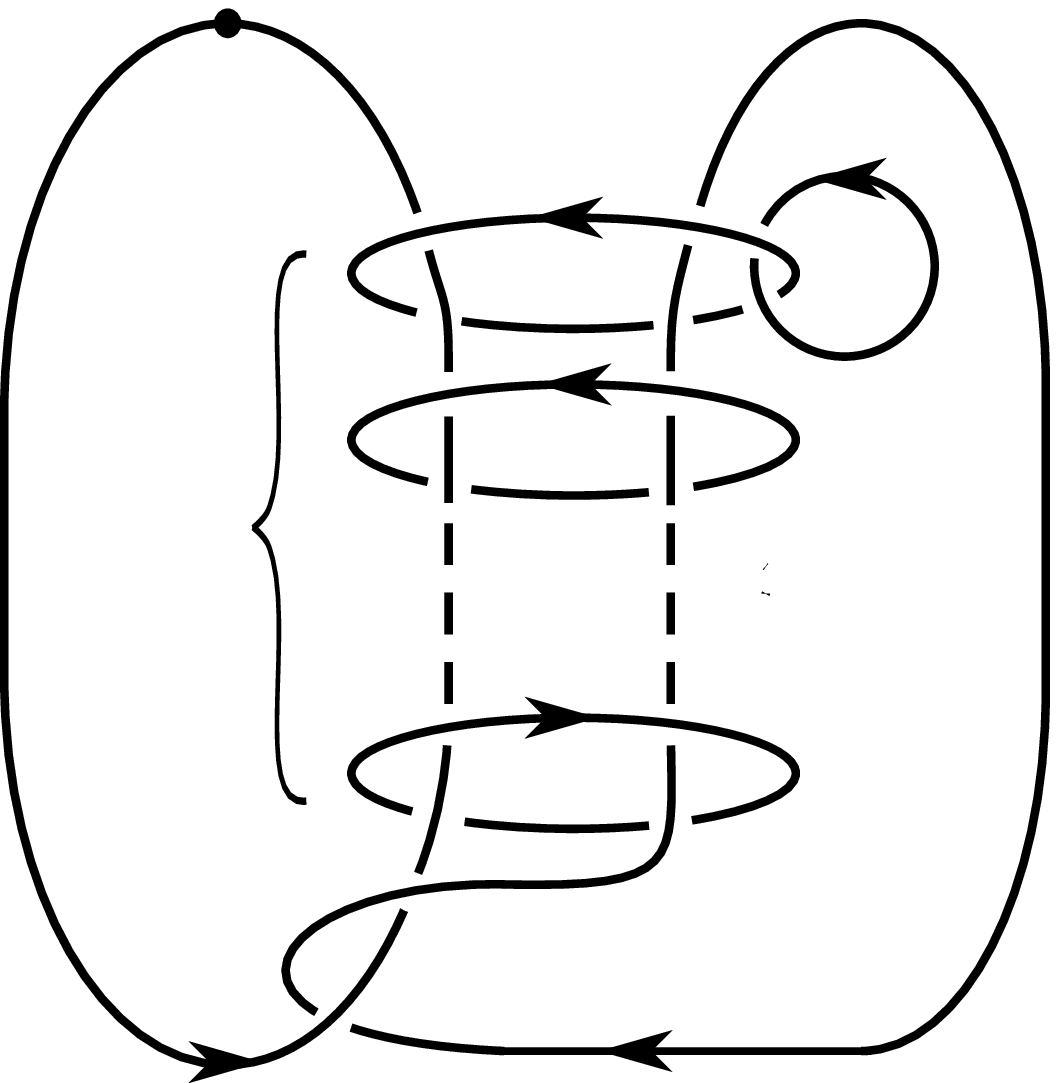}
\caption{ }
\labellist
\small\hair 2pt
\pinlabel $n-1$ at -330 180
\pinlabel $-n+1$ at -120 330
\pinlabel $-1$ at -200 280
\pinlabel $-1$ at -110 220
\pinlabel $-1$ at -140 190
\pinlabel $-1$ at -140 95
\pinlabel $n-2$ at -25 202
\pinlabel $-1$ at 135 267
\pinlabel $-1$ at 125 215
\pinlabel $-1$ at 122 115
\pinlabel $-1$ at 230 310
\pinlabel $-1$ at 320 260
\pinlabel $n-1$ at 280 175
\pinlabel $\text{twists}$ at 280 159
\endlabellist
\label{f:KirbyThmA7}
\end{minipage}
\hspace{0.3cm}
\begin{minipage}[b]{0.5\linewidth}
\centering
\includegraphics[height=65mm, width=60mm]{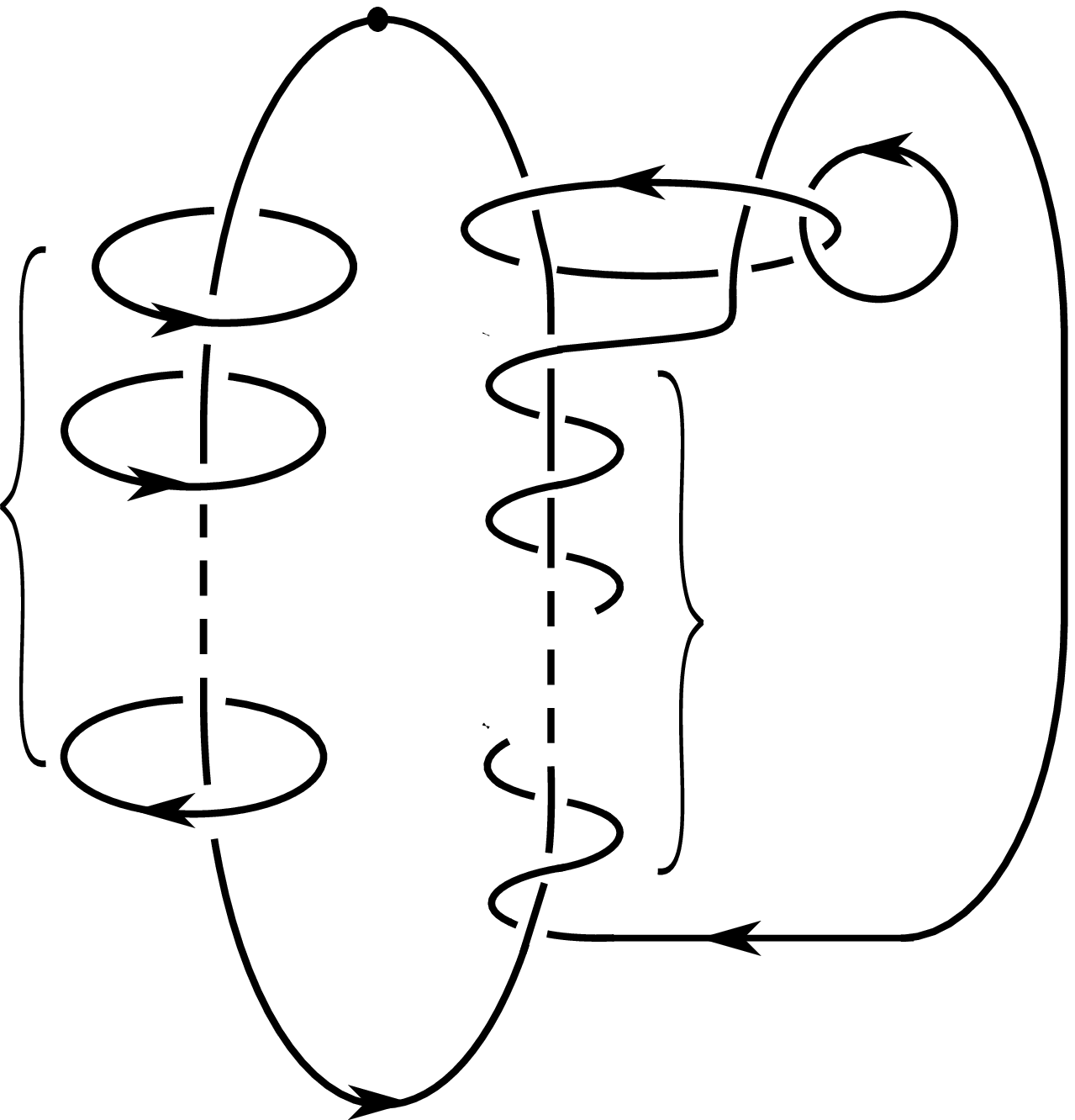}
\caption{ }
\label{f:KirbyThmA8}
\end{minipage}
\end{figure}

\begin{figure}[ht]
\begin{minipage}[b]{0.55\linewidth}
\centering
\includegraphics[scale=0.47]{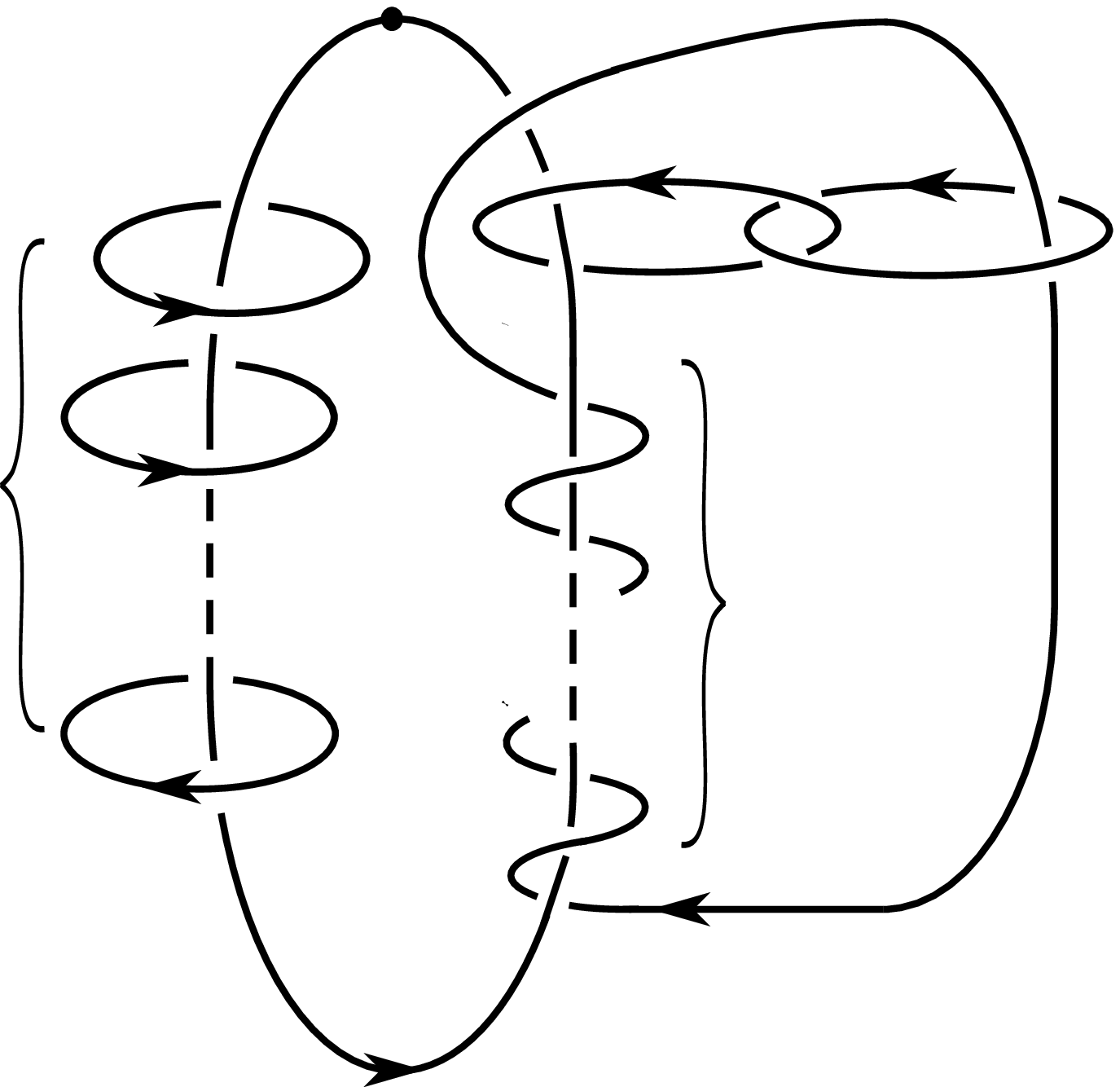}
\caption{ }
\labellist
\small\hair 2pt
\pinlabel $n-2$ at -490 175
\pinlabel $-1$ at -340 232
\pinlabel $-1$ at -345 180
\pinlabel $-1$ at -345 70
\pinlabel $-1$ at -210 290
\pinlabel $-1$ at -150 290
\pinlabel $0$ at -210 345
\pinlabel $n$ at -175 140
\pinlabel $\text{twists}$ at -175 124
\pinlabel $0$ at 190 330
\pinlabel $n-1$ at 55 380
\pinlabel $n$ at 190 177
\pinlabel $\text{twists}$ at 190 161
\endlabellist
\label{f:KirbyThmA9}
\end{minipage}
\hspace{0.3cm}
\begin{minipage}[b]{0.4\linewidth}
\centering
\includegraphics[scale=0.47]{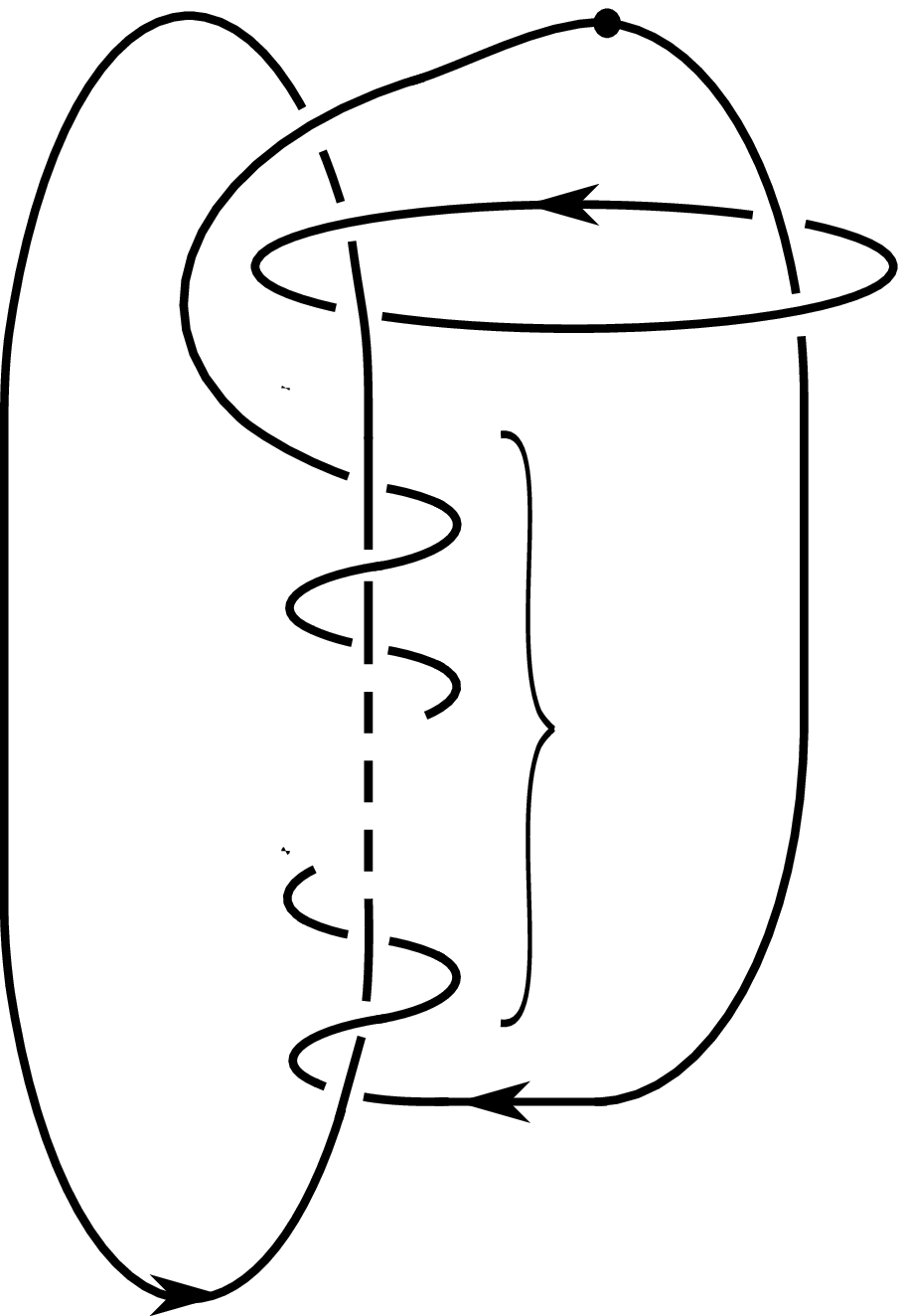}
\caption{$B_n$ with a two-handle}
\label{f:KirbyThmA10}
\end{minipage}
\end{figure}

\begin{figure}[ht]
\begin{minipage}[b]{0.36\linewidth}
\centering
\includegraphics[scale=0.42]{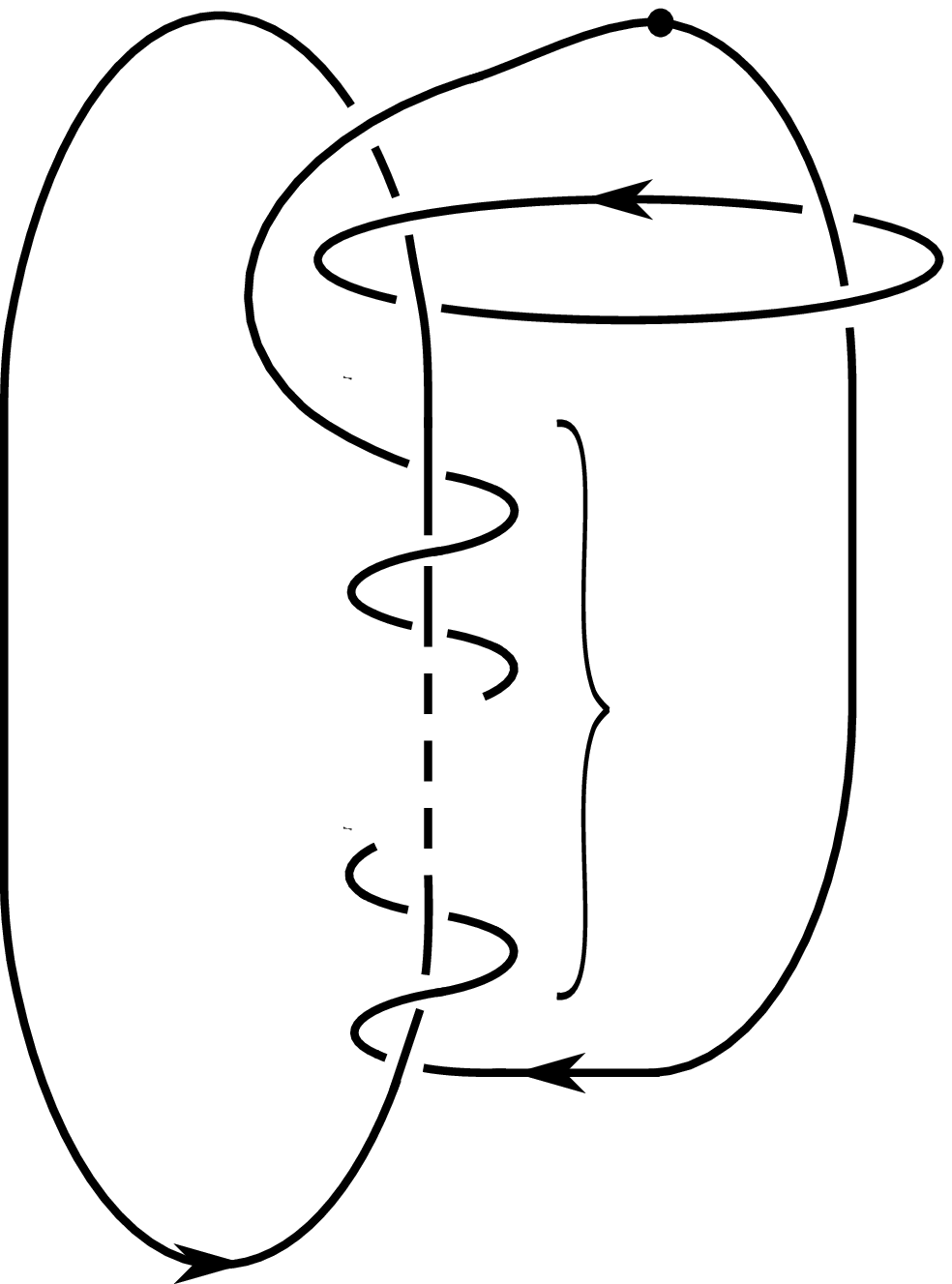}
\caption{ }
\labellist
\small\hair 2pt
\pinlabel $0$ at -150 190
\pinlabel $n-3$ at -370 218
\pinlabel $n-1$ at -150 105
\pinlabel $\text{twists}$ at -150 95
\pinlabel $0$ at 200 104
\pinlabel $-n-1$ at 330 104
\pinlabel $-n-1$ at 630 104
\endlabellist
\label{f:KirbyThmA11}
\end{minipage}
\begin{minipage}[b]{0.31\linewidth}
\centering
\includegraphics[height=40mm, width=35mm]{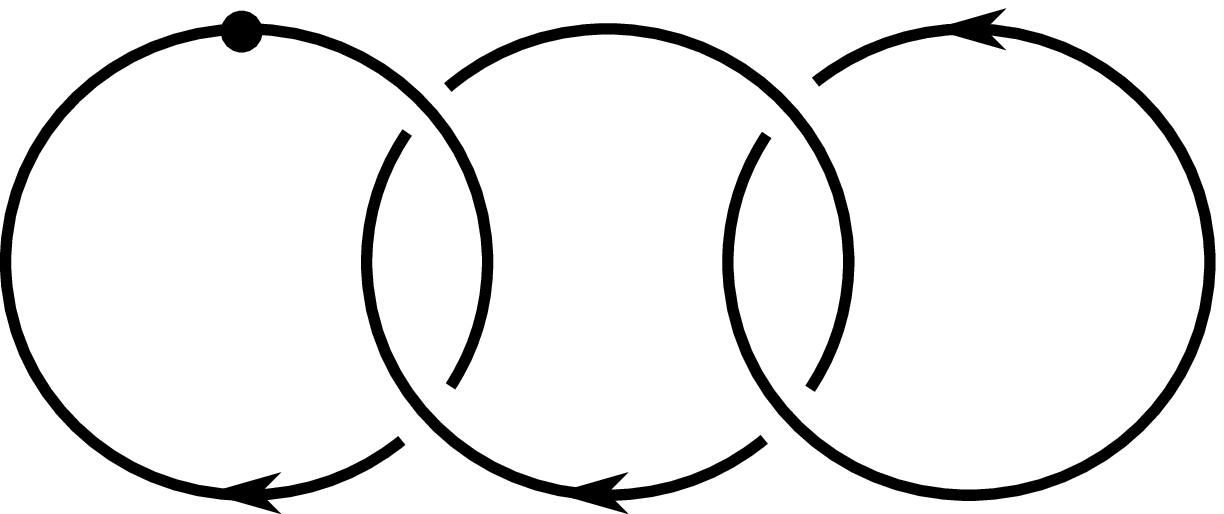}
\caption{ }
\label{f:KirbyThmA12}
\end{minipage}
\begin{minipage}[b]{0.31\linewidth}
\centering
\includegraphics[scale=0.28]{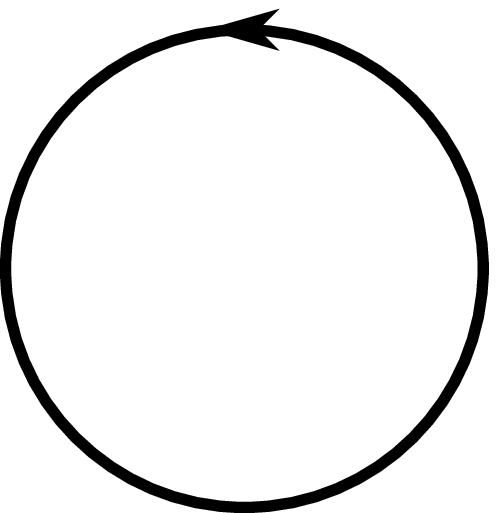}
\caption{$V_{-n-1}$}
\label{f:KirbyThmA13}
\end{minipage}
\end{figure}

Next, in Figure~\ref{f:KirbyThmA10} we perform the rational blow-down, thus replacing $C_n$ with $B_n$. This is done by first swapping the one-handle and the $0$-framed two-handle and then blowing down the $(n-1)$ spheres with self-intersection $(-1)$. Consequently, after rationally blowing down, we obtain $B_n$ with an additional $0$-framed two-handle. When we slide then $(n-1)$-framed two-handle of $B_n$ over that $0$-framed two-handle, we obtain Figure~\ref{f:KirbyThmA11}. We proceed to slide the same handle over the $0$-framed two-handle $(n-2)$ more times, and obtain Figure~\ref{f:KirbyThmA12}. Finally, we remove the cancelling $1/2$ handle pair, and obtain a single $(-n-1)$-framed two-handle, Figure~\ref{f:KirbyThmA13}, which is the manifold $V_{-n-1}$. Consequently, since to get from Figure~\ref{f:KirbyThmA10} to Figure~\ref{f:KirbyThmA13}, we only performed handleslides, it follows that $B_n \hookrightarrow V_{-n-1}$.    \end{proof}

\begin{cor}
\label{cor:smoothbn1}
For $n \geq 2$, the rational blow-up of $B_n \hookrightarrow V_{-n-1}$ is diffeomorphic to $V_{-n-1} \# (n-1)\overline{\C P^2}$.
\end{cor}

Corollary~\ref{cor:smoothbn1} follows directly from the proof of Theorem~\ref{thm:smoothbn1}. If we follow the Kirby moves backwards from Figure~\ref{f:KirbyThmA13} to Figure~\ref{f:KirbyThmA4}, it follows that if we start with a $V_{-n-1}$, and rationally blow up the $B_n \hookrightarrow V_{-n-1}$, then we end up with $V_{-n-1} \# (n-1)\overline{\C P^2}$.

\section{Proof of Theorem~\ref{thm:smoothbn2}}
\label{sec:pfthm2}

In this section we prove Theorem~\ref{thm:smoothbn2}.

\begin{proof}{\textit{Proof of Theorem~\ref{thm:smoothbn2}}.}

We prove Theorem~\ref{thm:smoothbn2} using similar Kirby calculus techniques as in the proof of Theorem~\ref{thm:smoothbn1}. (Note, the case $n=2$ is trivial and the case $n=3$ is covered in Theorem~\ref{thm:smoothbn1}, so here we can assume $n \geq 4$.) We start with the Kirby diagram for $V_{-4}$, Figure~\ref{f:KirbyThmB1}. We blow up $V_{-4}$ $(n-1)$ times, as seen in Figures~\ref{f:KirbyThmB2} through \ref{f:KirbyThmB5}, in such a manner that we end up with a plumbing tree of spheres as seen in Figure~\ref{f:KirbyThmB5}. This configuration of spheres is $C_n$ with an extra sphere $\Sigma'_{-1}$, with self-intersection $(-1)$, which intersects only with the first sphere with self-intersection $(-2)$, $S_2$, (compare with Figure~\ref{f:KirbyThmA4}).

\begin{figure}[ht]
\begin{minipage}[b]{0.32\linewidth}
\centering
\includegraphics[scale=0.28]{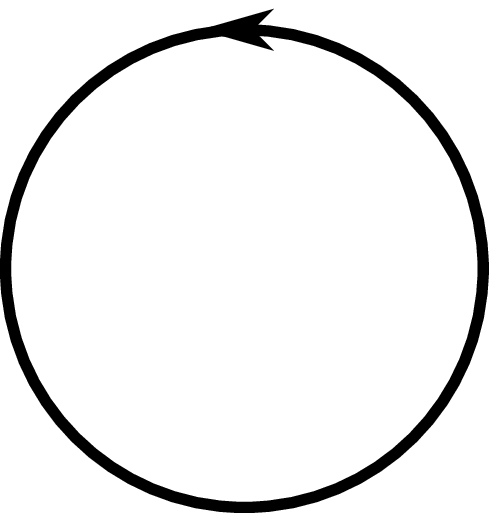}
\caption{$V_{-4}$}
\labellist
\small\hair 2pt
\pinlabel $-4$ at -260 130
\pinlabel $-5$ at 40 130
\pinlabel $-1$ at 200 130
\pinlabel $-6$ at 350 130
\pinlabel $-1$ at 500 130
\pinlabel $-2$ at 590 130
\endlabellist
\label{f:KirbyThmB1}
\end{minipage}
\begin{minipage}[b]{0.33\linewidth}
\centering
\includegraphics[height=20mm, width=30mm]{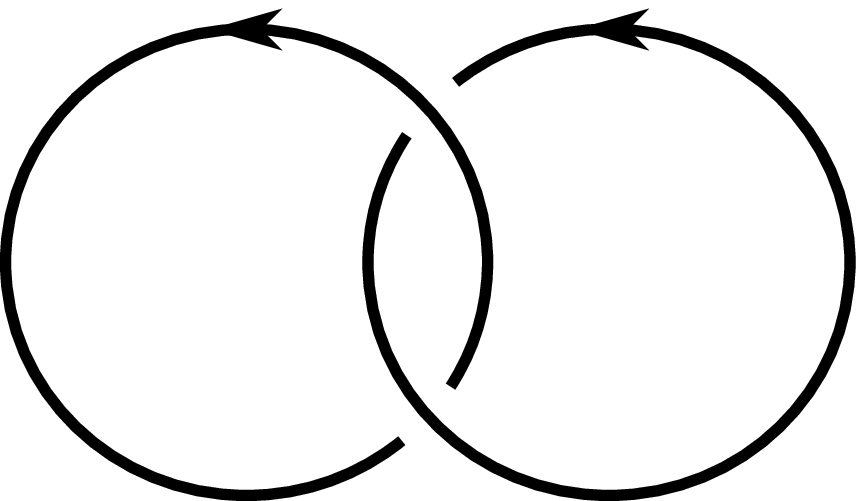}
\caption{$V_{-4} \# \overline{\C P^2}$}
\label{f:KirbyThmB2}
\end{minipage}
\begin{minipage}[b]{0.33\linewidth}
\centering
\includegraphics[scale=0.28]{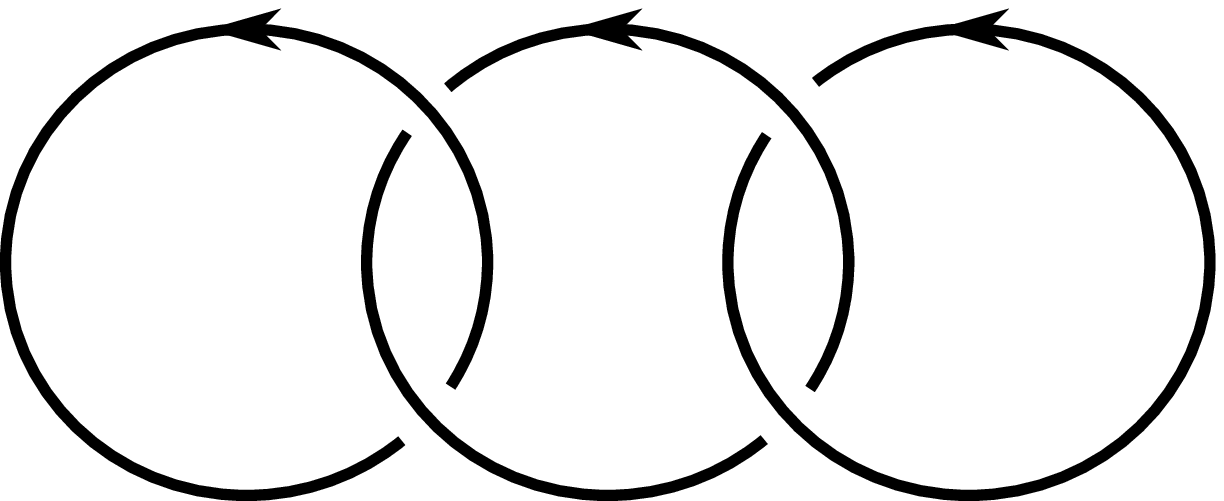}
\caption{$V_{-4} \# 2\overline{\C P^2}$}
\label{f:KirbyThmB3}
\end{minipage}
\end{figure}

\begin{figure}[ht!]
\labellist
\small\hair 2pt
\pinlabel $-n-2$ at 40 175
\pinlabel $-1$ at 230 170
\pinlabel $-2$ at 330 170
\pinlabel $-2$ at 430 170
\pinlabel $-2$ at 770 170
\pinlabel $n-3$ at 520 -5
\endlabellist
\centering
\includegraphics[height=25mm, width=90mm]{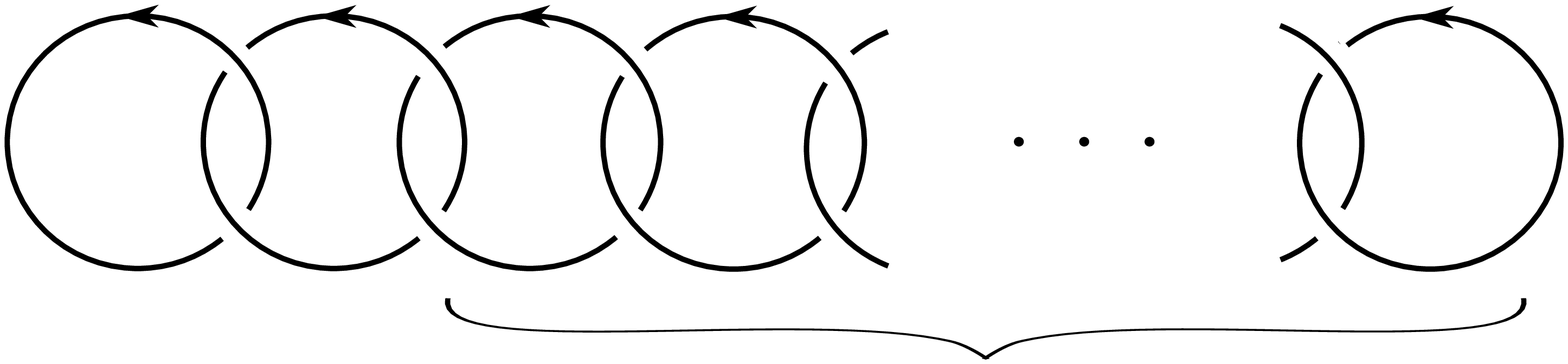}
\caption{$V_{-4} \# (n-2)\overline{\C P^2}$}
\label{f:KirbyThmB4}
\end{figure}

\begin{figure}[ht!]
\labellist
\small\hair 2pt
\pinlabel $-n-2$ at 20 265
\pinlabel $-2$ at 210 260
\pinlabel $-2$ at 310 260
\pinlabel $-2$ at 420 260
\pinlabel $-2$ at 750 260
\pinlabel $-1$ at 230 5
\pinlabel $n-2$ at 450 330
\endlabellist
\centering
\includegraphics[height=40mm, width=90mm]{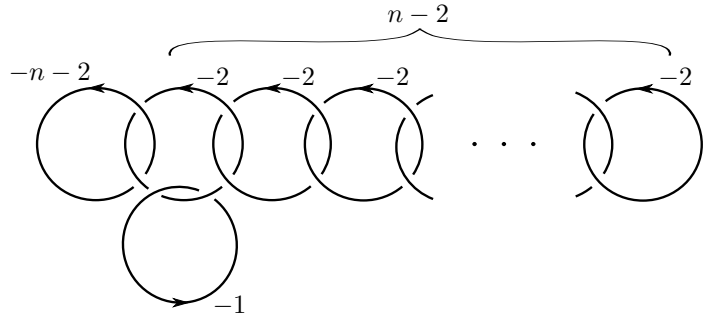}
\caption{$V_{-4} \# (n-1)\overline{\C P^2}$}
\label{f:KirbyThmB5}
\end{figure}

\begin{figure}[ht!]
\labellist
\small\hair 2pt
\pinlabel $-n-2$ at 220 35
\pinlabel $-1$ at 375 40
\pinlabel $-2$ at 325 170
\pinlabel $-2$ at 375 170
\pinlabel $-2$ at 530 170
\pinlabel $n-2$ at 420 210
\pinlabel $-1$ at 250 175
\endlabellist
\centering
\includegraphics[height=45mm, width=125mm]{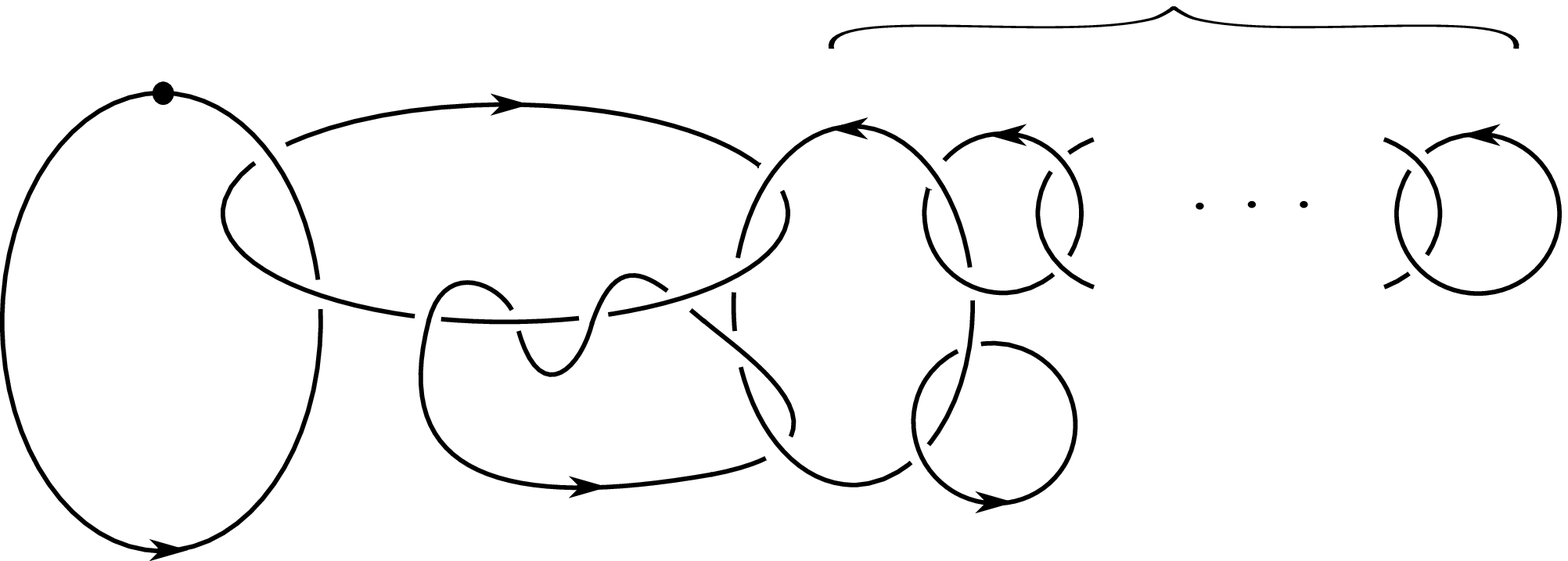}
\caption{ }
\label{f:KirbyThmB6}
\end{figure}

\begin{figure}[ht!]
\labellist
\small\hair 2pt
\pinlabel $-n+1$ at 500 60
\pinlabel $-2$ at 265 165
\pinlabel $-2$ at 310 165
\pinlabel $-2$ at 455 163
\pinlabel $-1$ at 280 60
\pinlabel $n-2$ at 350 195
\pinlabel $-1$ at 200 170
\endlabellist
\centering
\includegraphics[height=45mm, width=120mm]{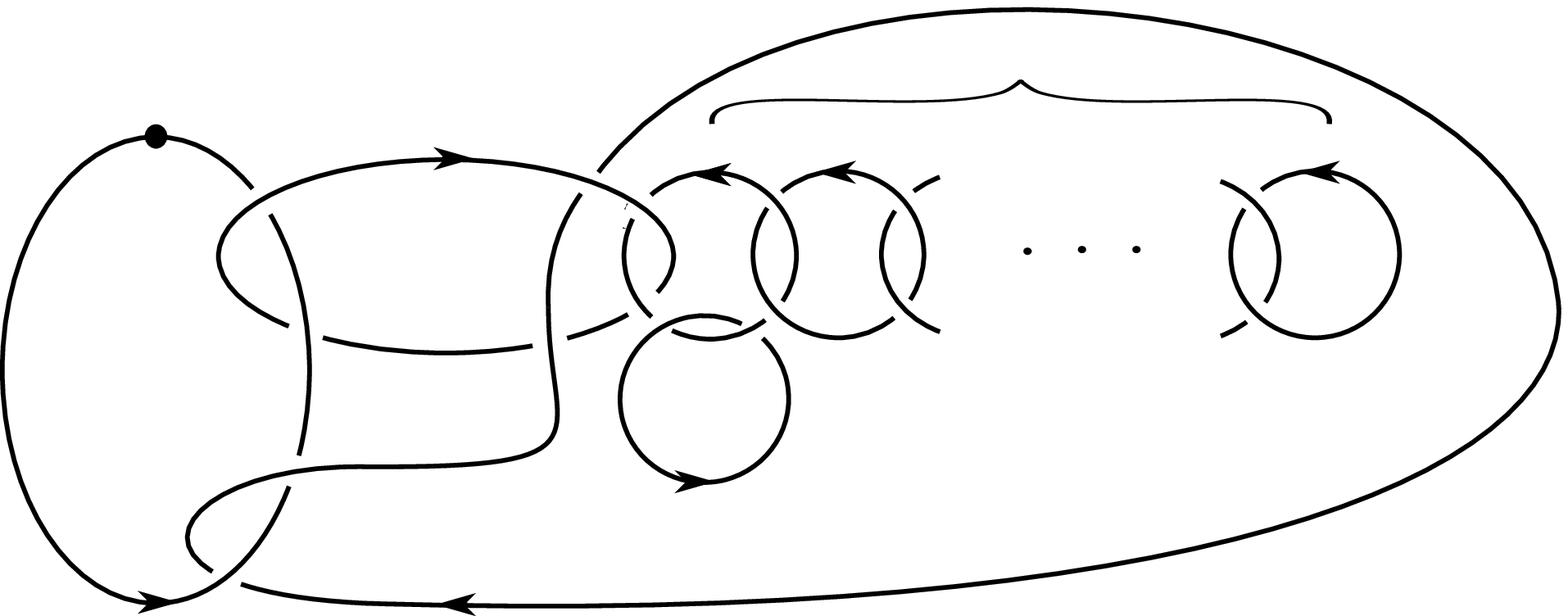}
\caption{ }
\label{f:KirbyThmB7}
\end{figure}

\begin{figure}[ht!]
\labellist
\small\hair 2pt
\pinlabel $-n+1$ at 500 60
\pinlabel $-1$ at 265 95
\pinlabel $-2$ at 310 183
\pinlabel $-2$ at 455 183
\pinlabel $-1$ at 220 77
\pinlabel $n-3$ at 365 90
\pinlabel $-1$ at 200 190
\endlabellist
\centering
\includegraphics[height=45mm, width=120mm]{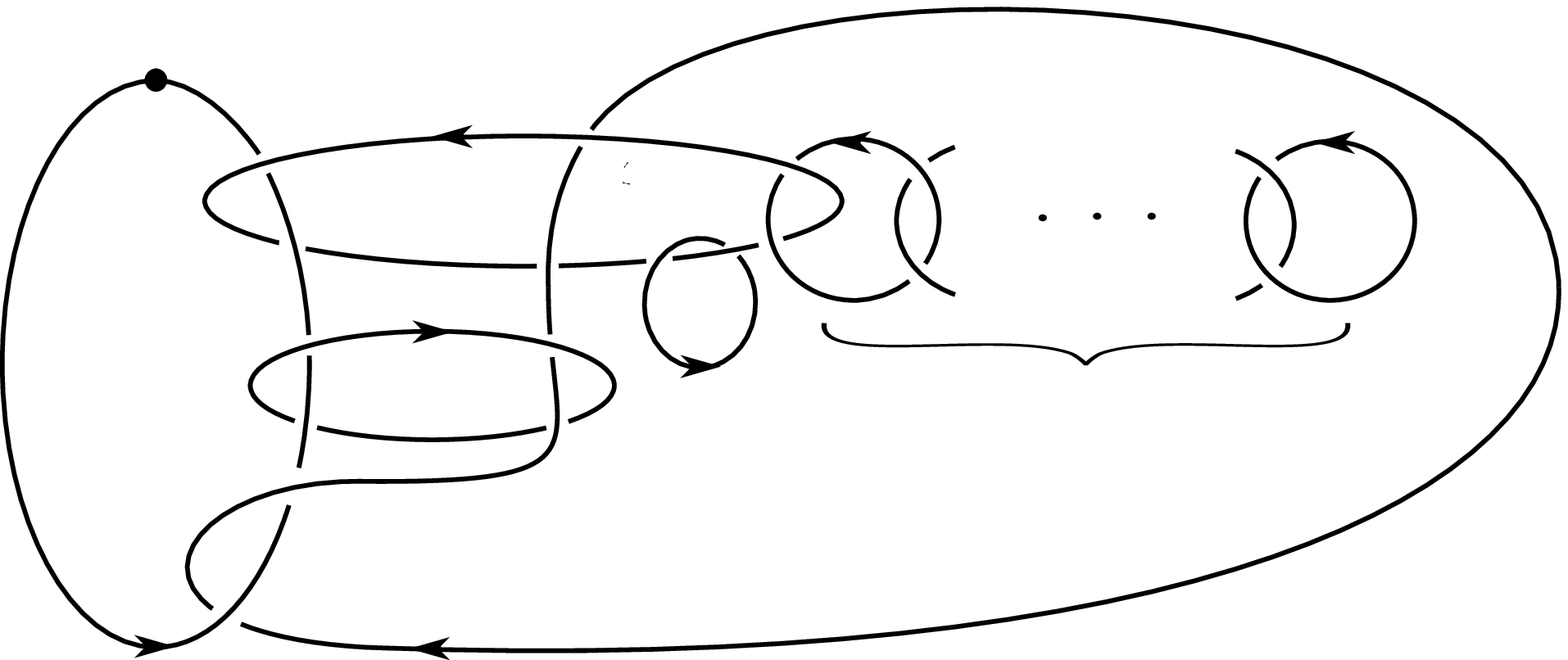}
\caption{ }
\label{f:KirbyThmB8}
\end{figure}

\begin{figure}[ht!]
\labellist
\small\hair 2pt
\pinlabel $-n+1$ at 480 60
\pinlabel $-2$ at 310 260
\pinlabel $-2$ at 455 260
\pinlabel $n-4$ at 372 160
\pinlabel $-1$ at 260 125
\pinlabel $-1$ at 90 100
\pinlabel $-1$ at 90 165
\pinlabel $-1$ at 90 230
\endlabellist
\centering
\includegraphics[height=55mm, width=120mm]{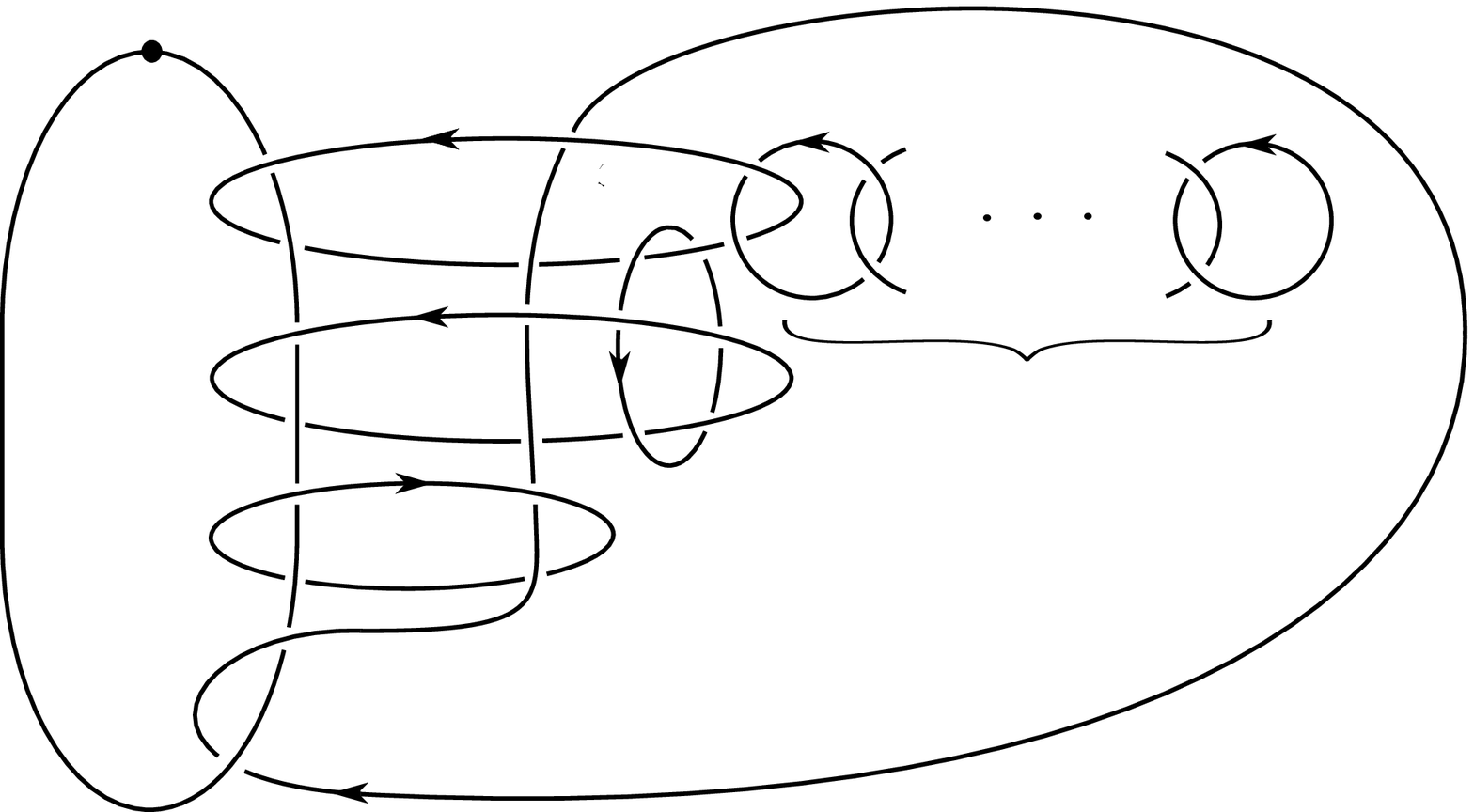}
\caption{ }
\label{f:KirbyThmB9}
\end{figure}

\begin{figure}[ht]
\begin{minipage}[b]{0.48\linewidth}
\centering
\includegraphics[scale=0.42]{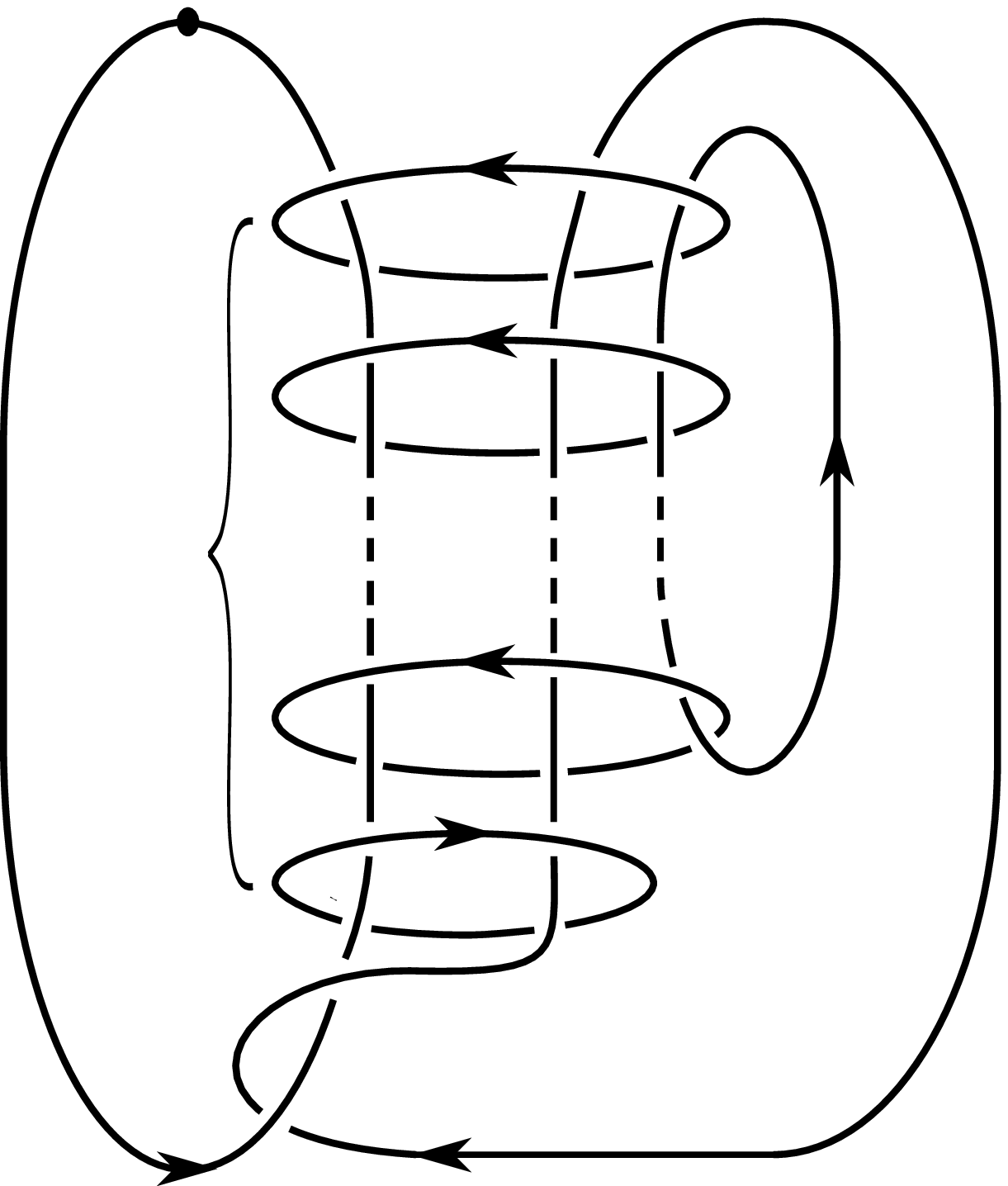}
\caption{ }
\labellist
\small\hair 2pt
\pinlabel $n-1$ at -375 235
\pinlabel $-n+1$ at -35 330
\pinlabel $-1$ at -130 160
\pinlabel $-1$ at -280 380
\pinlabel $-1$ at -320 315
\pinlabel $-1$ at -320 200
\pinlabel $-1$ at -320 143
\pinlabel $n-2$ at 90 250
\pinlabel $-1$ at 230 370
\pinlabel $-1$ at 150 300
\pinlabel $-1$ at 150 190
\pinlabel $-1$ at 110 130
\pinlabel $-1$ at 320 140
\pinlabel $-n+2$ at 330 45
\endlabellist
\label{f:KirbyThmB10}
\end{minipage}
\hspace{0.3cm}
\begin{minipage}[b]{0.47\linewidth}
\centering
\includegraphics[scale=0.43]{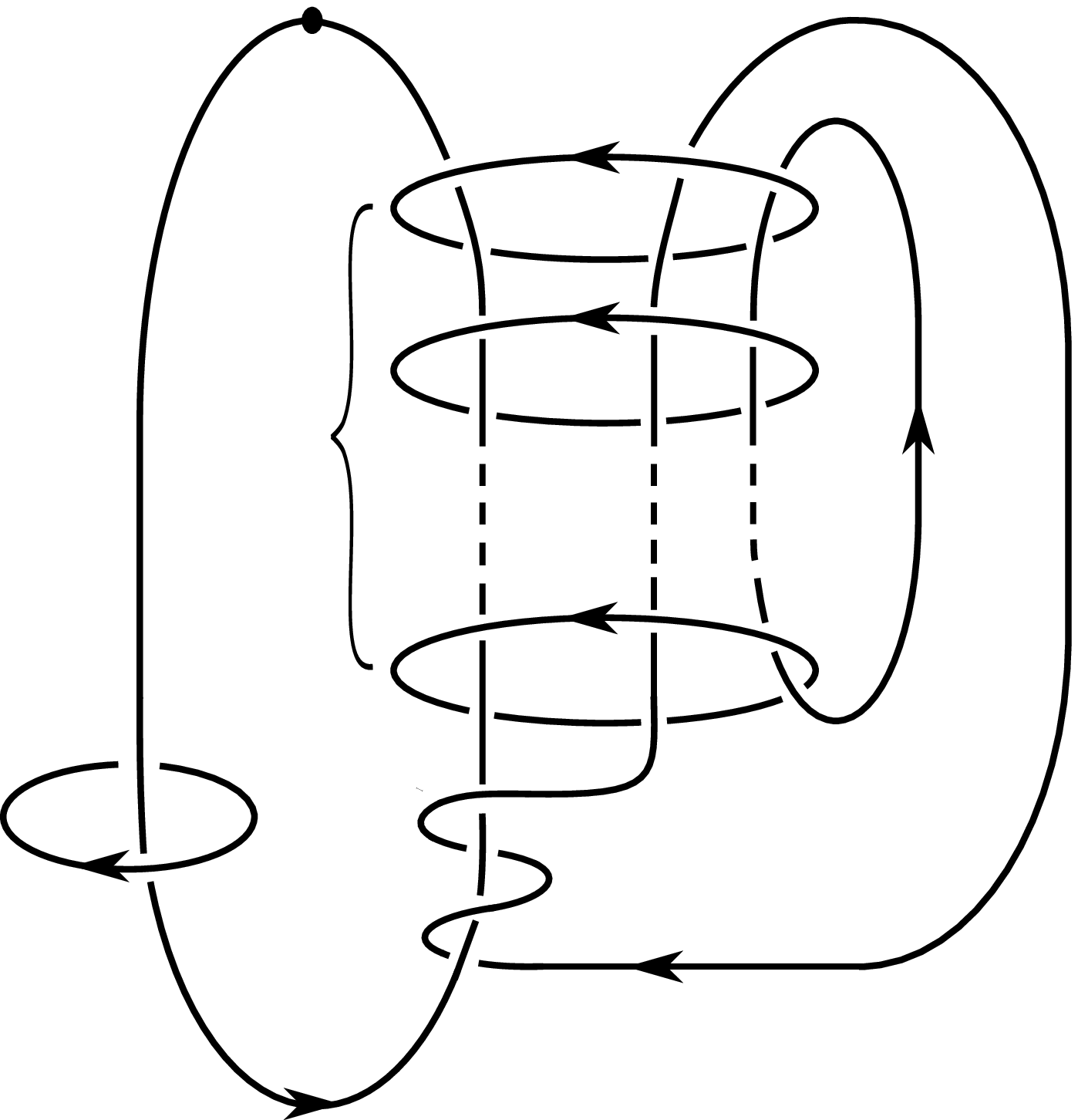}
\caption{ }
\label{f:KirbyThmB11}
\end{minipage}
\end{figure}

As was done in the proof of the previous theorem, we proceed with a series of Kirby moves that will change the presentation of $C_n$ in Figure~\ref{f:KirbyThmB5}, from a linear plumbing of spheres, Figure~\ref{f:KirbyCnCircles1}, to the one in Figure~\ref{f:KirbyCnCircles2}. We start by adding a cancelling $1/2$-handle pair in Figure~\ref{f:KirbyThmB6}. We proceed by sliding the $(-n-2)$-framed two-handle over the two-handle which was added in the previous step (Figure~\ref{f:KirbyThmB7}). Following this, we perform $(n-2)$ handleslides in order to slide off the $(-2)$-framed two-handles in Figures~\ref{f:KirbyThmB8}-\ref{f:KirbyThmB10}. As a result, the $(-1)$-framed two-handle corresponding to the sphere $\Sigma'_{-1}$, intersects once with each of the spheres corresponding to the $(n-2)$ $(-1)$-framed two-handles, as seen in Figure~\ref{f:KirbyThmB10}. Next, we slide the $(-n+1)$-framed two-handle off of each of the $(n-1)$ $(-1)$-framed two-handles, Figures~\ref{f:KirbyThmB11} and \ref{f:KirbyThmB12}. Consequently, in Figure~\ref{f:KirbyThmB12} we obtain a presentation of $C_n$ as in Figure~\ref{f:KirbyCnCircles2}, with the extra sphere $\Sigma'_{-1}$.

\begin{figure}[ht]
\begin{minipage}[b]{0.51\linewidth}
\centering
\includegraphics[scale=0.37]{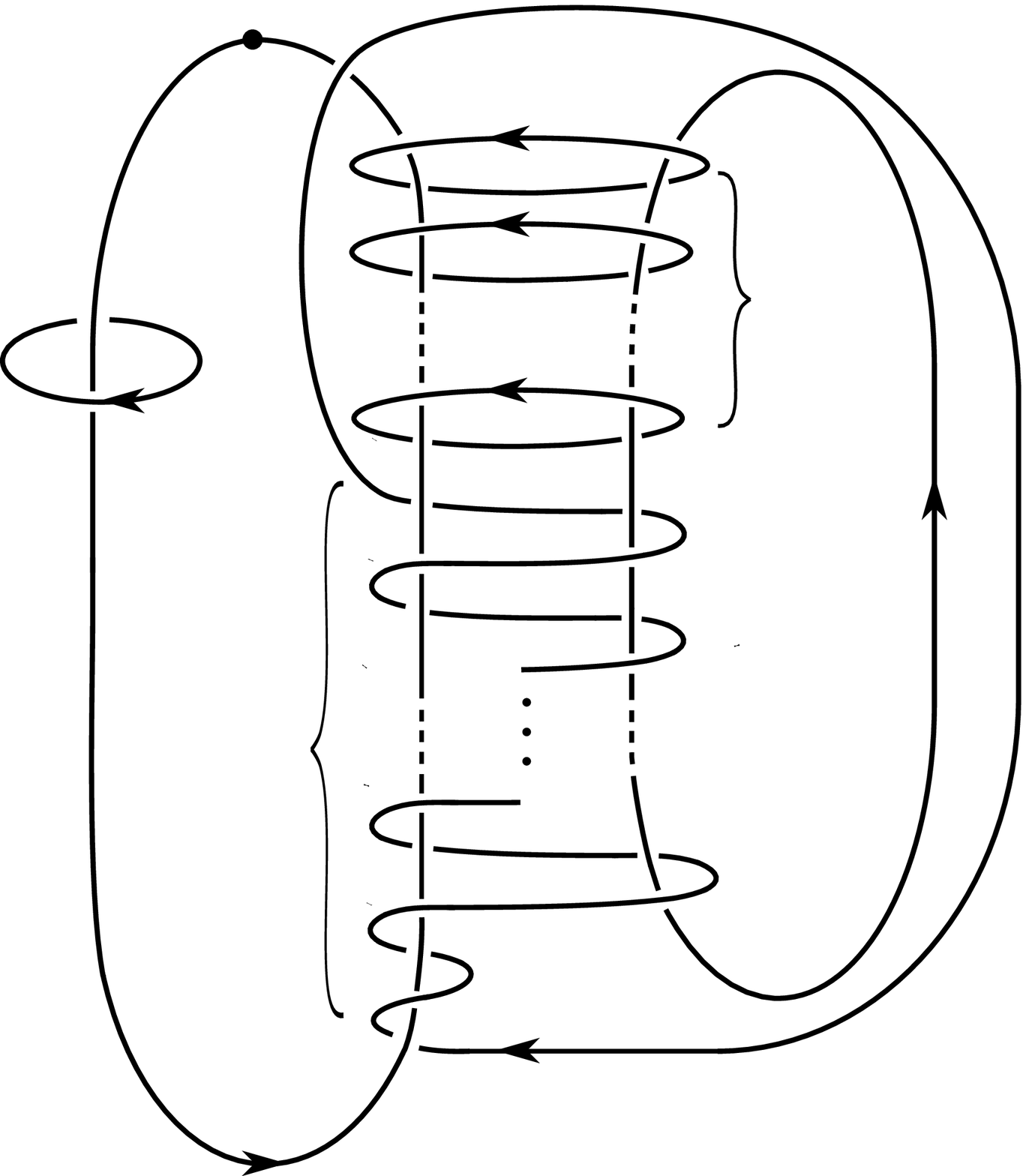}
\caption{ }
\labellist
\small\hair 2pt
\pinlabel $n-2$ at -148 410
\pinlabel $-1$ at -260 500
\pinlabel $-1$ at -215 450
\pinlabel $-1$ at -225 370
\pinlabel $-1$ at -470 350
\pinlabel $-1$ at -240 80
\pinlabel $0$ at -115 525
\pinlabel $n$ at -450 197
\pinlabel $\text{twists}$ at -450 179
\pinlabel $n-3$ at 240 100
\pinlabel $n-1$ at 80 550
\pinlabel $n$ at 70 200
\pinlabel $\text{twists}$ at 70 184
\pinlabel $n-2$ at 160 432
\pinlabel $\text{twists}$ at 160 416
\endlabellist
\label{f:KirbyThmB12}
\end{minipage}
\hspace{0.3cm}
\begin{minipage}[b]{0.44\linewidth}
\centering
\includegraphics[scale=0.37]{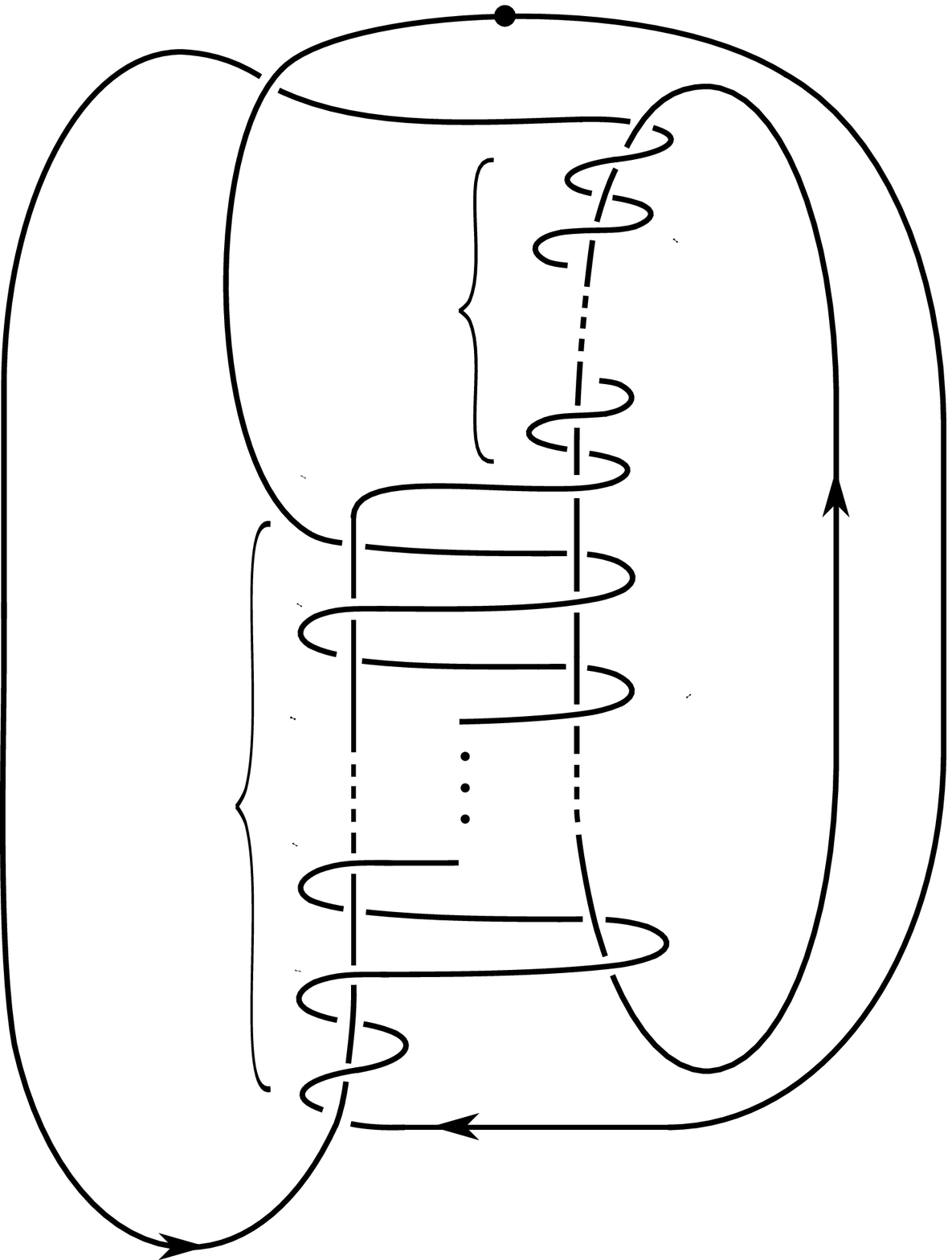}
\caption{$B_n$ with a two-handle}
\label{f:KirbyThmB13}
\end{minipage}
\end{figure}

\begin{figure}[ht]
\begin{minipage}[b]{0.48\linewidth}
\centering
\includegraphics[scale=0.5]{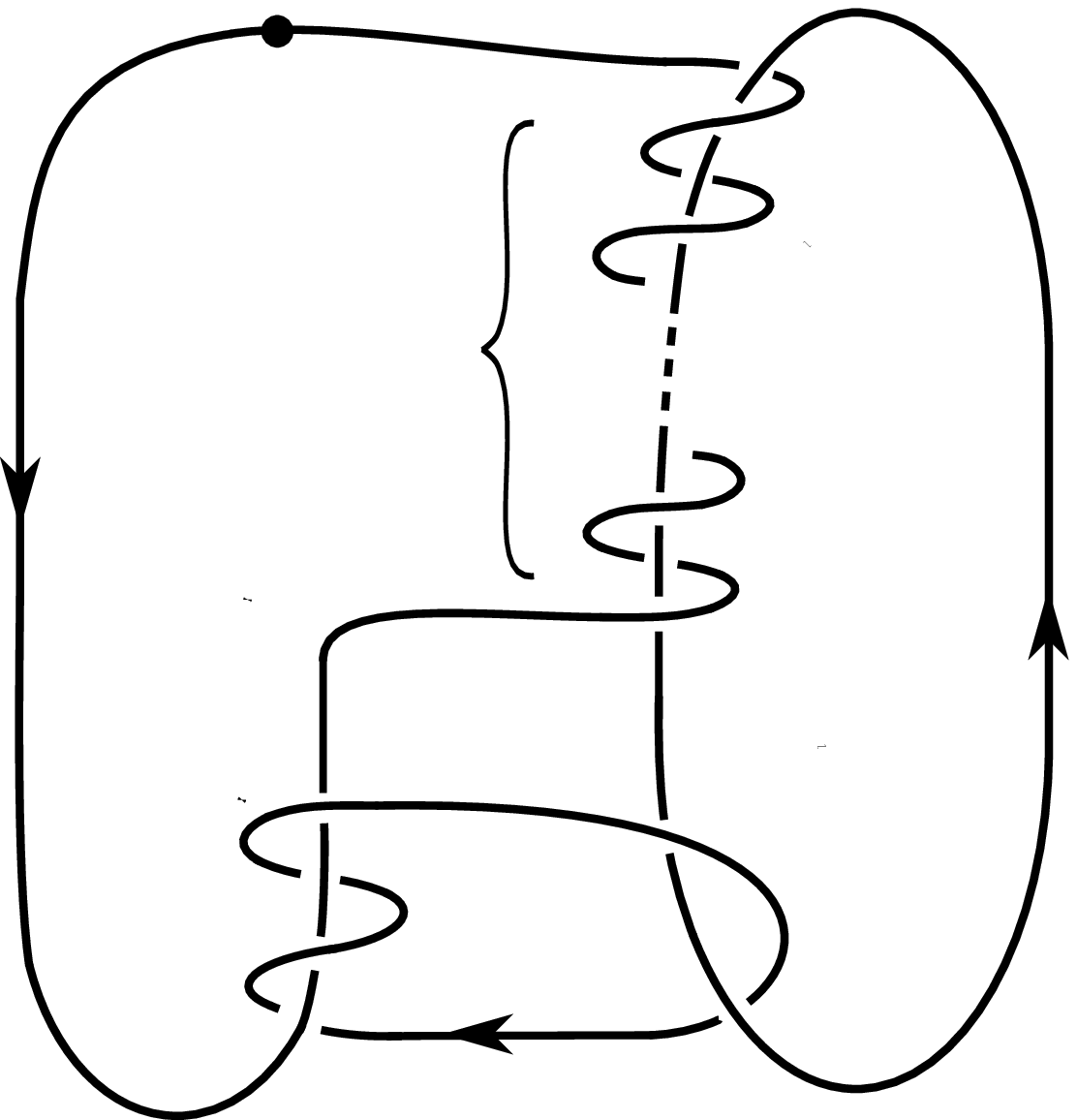}
\caption{ }
\labellist
\small\hair 2pt
\pinlabel $n-3$ at -115 342
\pinlabel $0$ at -210 15
\pinlabel $n-2$ at -260 242
\pinlabel $\text{twists}$ at -260 225
\pinlabel $n-5$ at 255 332
\pinlabel $0$ at 160 15
\pinlabel $n-4$ at 110 235
\pinlabel $\text{twists}$ at 110 218
\endlabellist
\label{f:KirbyThmB14}
\end{minipage}
\hspace{0.3cm}
\begin{minipage}[b]{0.47\linewidth}
\centering
\includegraphics[scale=0.5]{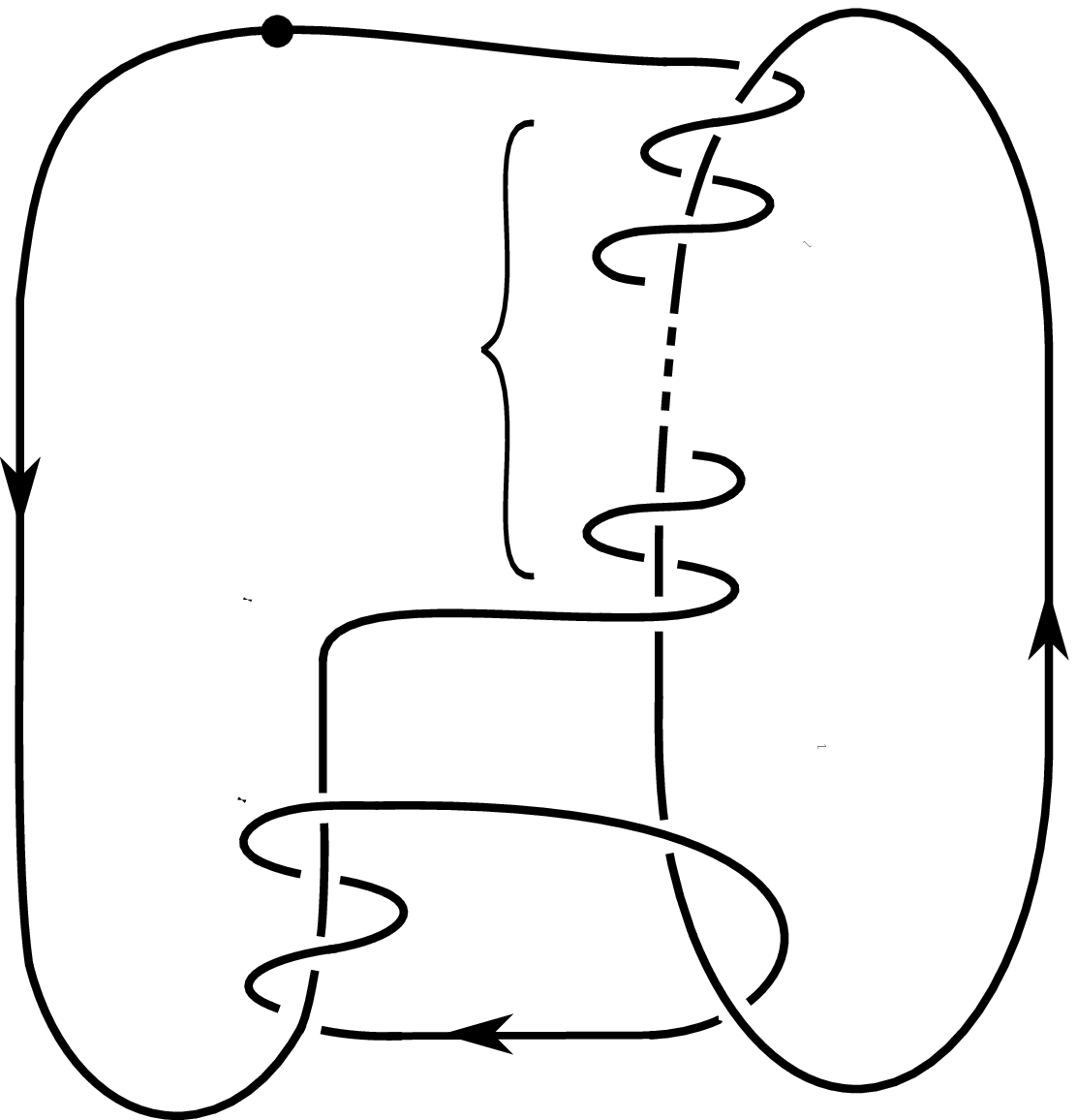}
\caption{ }
\label{f:KirbyThmB15}
\end{minipage}
\end{figure}

Next, we perform the rational blow-down procedure, by exchanging the one-handle and the $0$-framed two-handle, and blowing down along the $(n-1)$ spheres with self-intersection $(-1)$, and obtain the Kirby diagram of $B_n$ with an additional $(n-3)$-framed two-handle, Figure~\ref{f:KirbyThmB13}. Next, we slide the $(n-1)$-framed two-handle over the $(n-3)$-framed two-handle and obtain the Kirby diagram in Figure~\ref{f:KirbyThmB14}, with the $(n-1)$-framed two-handle becoming a 0-framed two-handle. At this point, the unknot corresponding to the $(n-3)$-framed two-handle is linked with the unknot corresponding to the one-handle with $(n-2)$ twists. If we slide off that $(n-3)$-framed two-handle off of the the $0$-framed two-handle, then we knock down the amount of twists that the unknot corresponding to the $(n-3)$-framed two-handle is linked with the unknot corresponding to the one-handle by $2$, thus obtaining Figure~\ref{f:KirbyThmB15}. 

\begin{figure}[ht]
\begin{minipage}[b]{0.4\linewidth}
\centering
\includegraphics[scale=0.55]{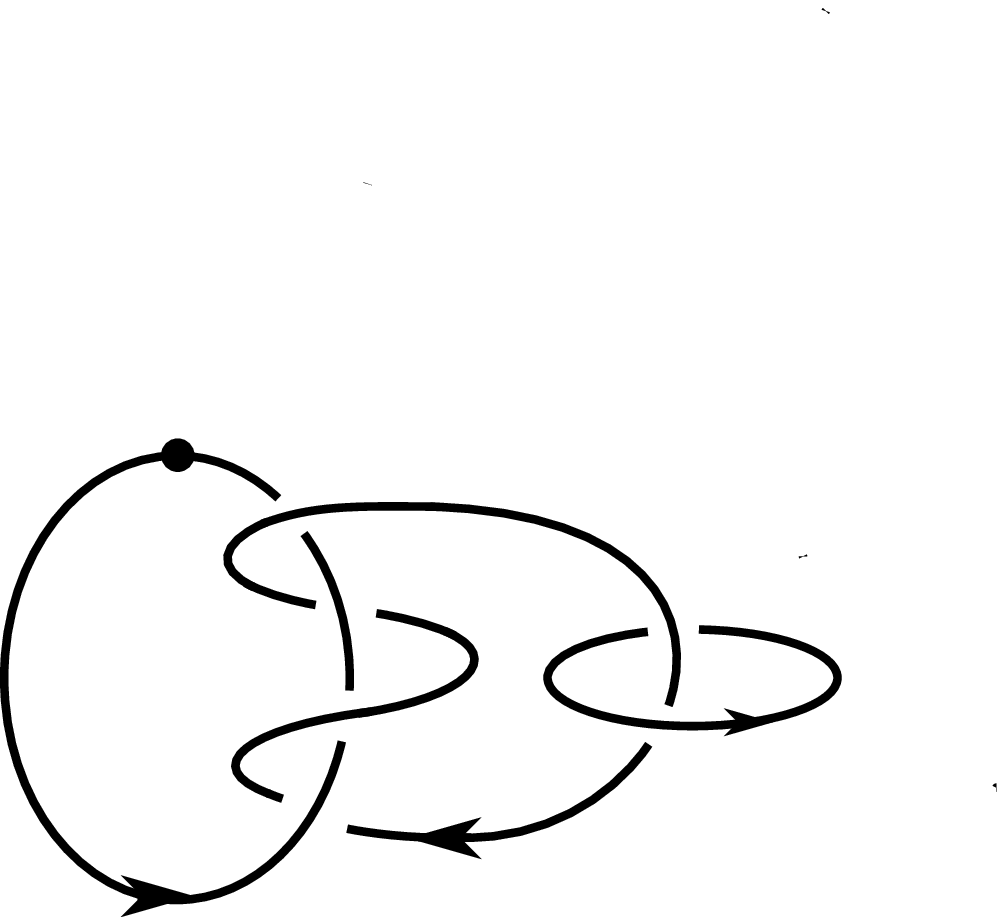}
\caption{ }
\labellist
\small\hair 2pt
\pinlabel $-1$ at -80 45
\pinlabel $0$ at -140 15
\pinlabel $-1$ at 300 58
\pinlabel $0$ at 70 153
\pinlabel $2$ at 122 75
\endlabellist
\label{f:KirbyThmB16even}
\end{minipage}
\hspace{0.3cm}
\begin{minipage}[b]{0.55\linewidth}
\centering
\includegraphics[scale=0.55]{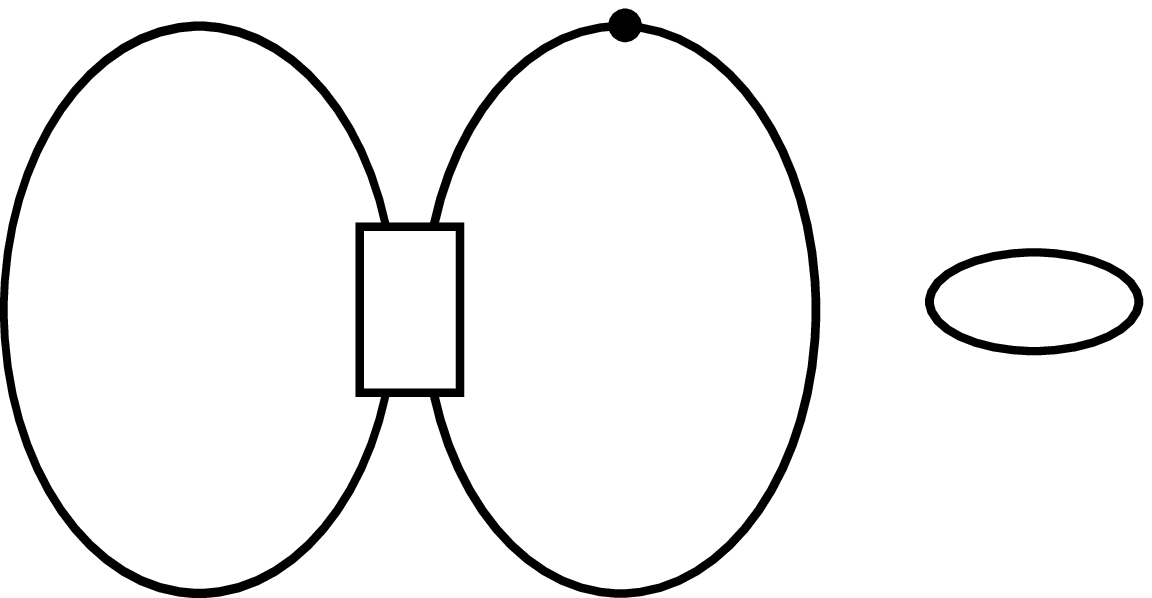}
\caption{$B_2 \# \overline{\C P^2}$}
\label{f:KirbyThmB17even}
\end{minipage}
\end{figure}

If $n$ is even, then after $\frac{n-2}{2}$ such handleslides we obtain the diagram in Figure~\ref{f:KirbyThmB16even}, (equivalent to the one in Figure~\ref{f:KirbyThmB17even}), which is just $B_2$ blown up once, i.e. $B_2 \# \overline{\C P^2}$. Consequently, if we start with $B_2 \# \overline{\C P^2}$, and follow the Kirby moves backwards from Figure~\ref{f:KirbyThmB17even} to Figure~\ref{f:KirbyThmB13}, then we see that $B_n \hookrightarrow B_2 \# \overline{\C P^2}$, for $n$ even.

\begin{figure}[ht]
\begin{minipage}[b]{0.48\linewidth}
\centering
\includegraphics[scale=0.45]{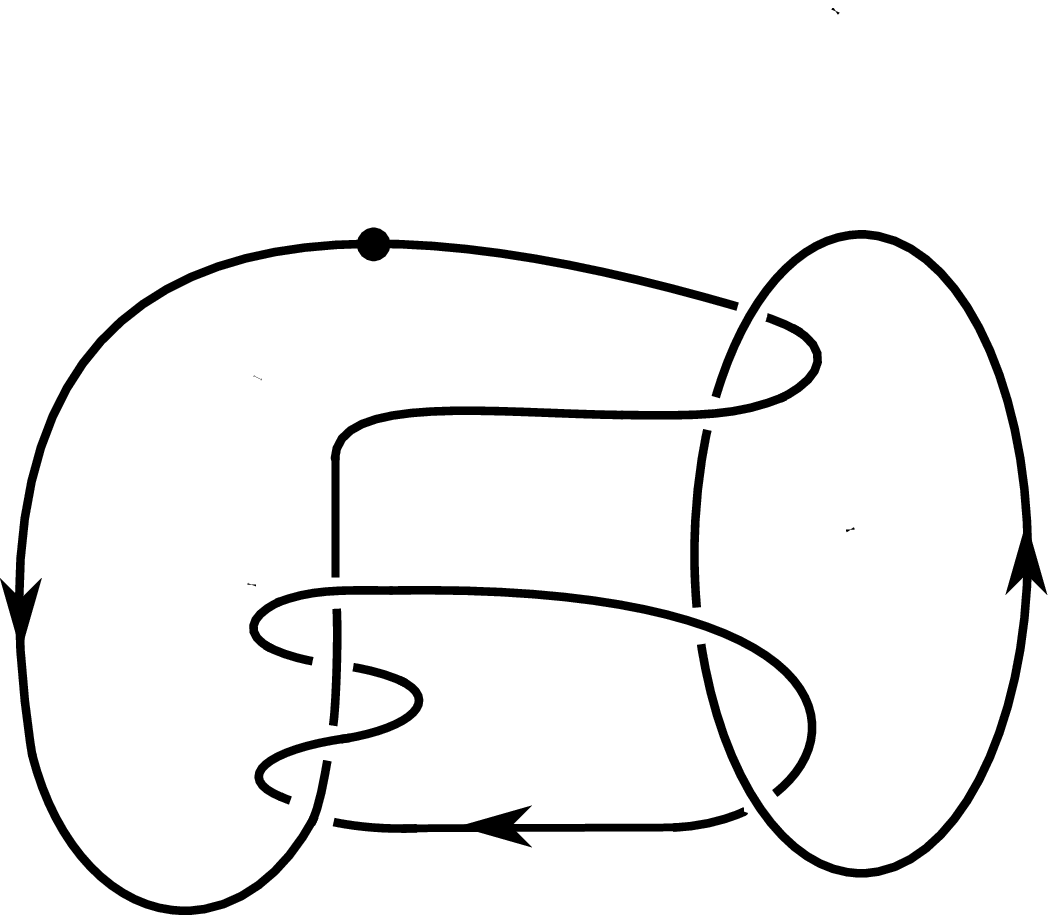}
\caption{ }
\labellist
\small\hair 2pt
\pinlabel $0$ at -160 253
\pinlabel $0$ at -240 15
\pinlabel $0$ at 235 258
\pinlabel $-2$ at 160 50
\endlabellist
\label{f:KirbyThmB16odd}
\end{minipage}
\hspace{0.3cm}
\begin{minipage}[b]{0.47\linewidth}
\centering
\includegraphics[height=35mm, width=50mm]{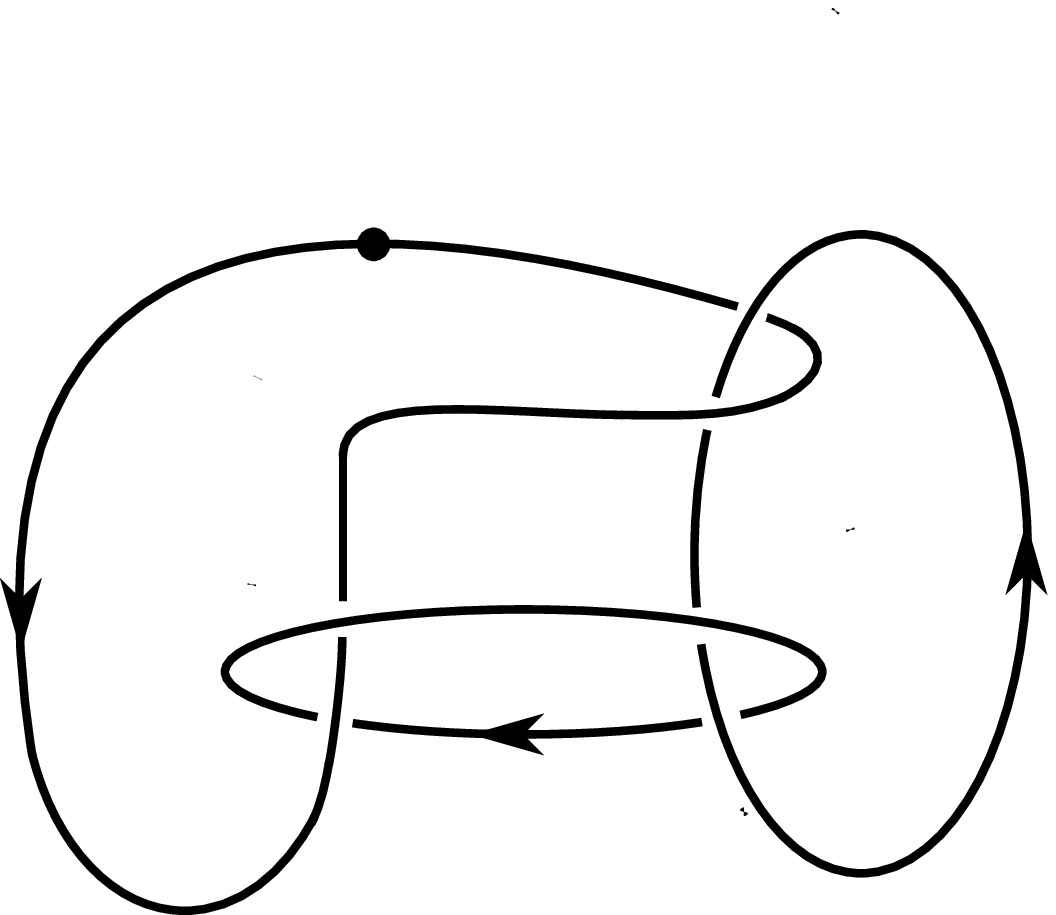}
\caption{ }
\label{f:KirbyThmB17odd}
\end{minipage}
\end{figure}

\begin{figure}[ht]
\begin{minipage}[b]{0.48\linewidth}
\centering
\includegraphics[scale=0.36]{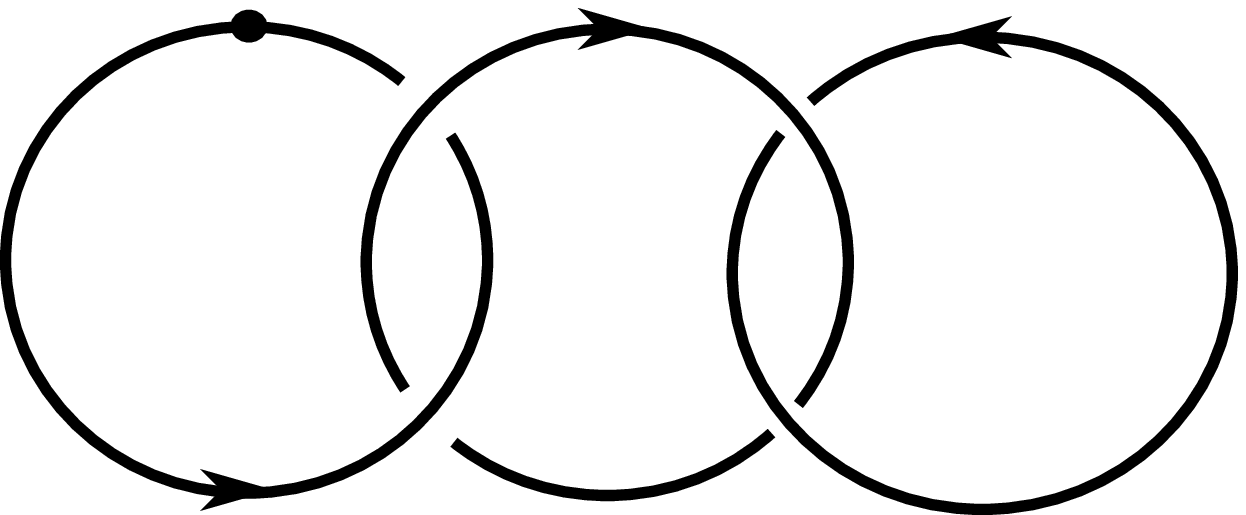}
\caption{ }
\labellist
\small\hair 2pt
\pinlabel $0$ at -280 135
\pinlabel $-4$ at -205 135
\pinlabel $-4$ at 105 130
\endlabellist
\label{f:KirbyThmB18odd}
\end{minipage}
\hspace{0.3cm}
\begin{minipage}[b]{0.47\linewidth}
\centering
\includegraphics[height=25mm, width=25mm]{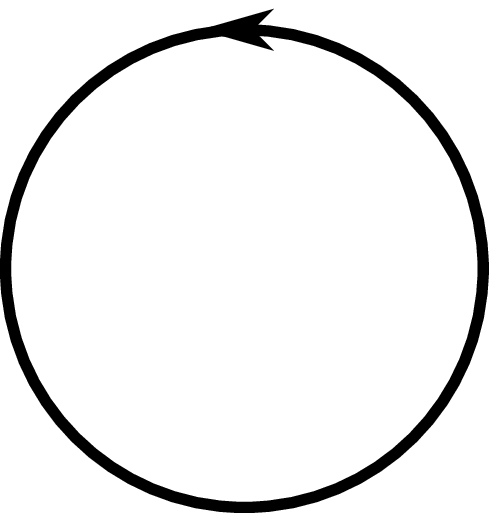}
\caption{$V_{-4}$}
\label{f:KirbyThmB19odd}
\end{minipage}
\end{figure}

If $n$ is odd, then if we start with the diagram in Figure~\ref{f:KirbyThmB14}, and slide off the $(n-3)$-framed two-handle $\frac{n-3}{2}$ times, we obtain the diagram in Figure~\ref{f:KirbyThmB16odd}. Following this, we slide the $0$-framed two handle (the one on the bottom of the diagram), over the other $0$-framed two-handle and obtain the diagram in Figure~\ref{f:KirbyThmB17odd}. We then perform another handleslide, and slide off the $(-2)$-framed two-handle off of the $0$-framed two-handle and get the diagram in Figure~\ref{f:KirbyThmB18odd}. Finally, we remove the cancelling $1/2$-handle pair and are left with one $(-4)$-framed two-handle, Figure~\ref{f:KirbyThmB19odd}, which represents the manifold $V_{-4}$. Consequently, if we follow the Kirby moves backwards from Figure~\ref{f:KirbyThmB19odd} to Figure~\ref{f:KirbyThmB13} (skipping Figures~\ref{f:KirbyThmB16even} and \ref{f:KirbyThmB17even}, as these are for the case when $n$ is even), then we can conclude that $B_n \hookrightarrow V_{-4}$ for $n$ odd.    \end{proof}

The difference between the embeddings in Theorem~\ref{thm:smoothbn2} with $n$ odd and even occurs because for $n$ odd the rational homology balls $B_n$ are spin and for $n$ even they are not.

\begin{cor}
\label{cor:smoothbn2}
For odd $n \geq 3$, the rational blow-up of $B_n \hookrightarrow V_{-4}$ is diffeomorphic to $V_{-4} \# (n-1)\overline{\C P^2}$.
\end{cor}

Similarly to the proof of Corollary~\ref{cor:smoothbn1}, Corollary~\ref{cor:smoothbn2} follows directly from the proof of Theorem~\ref{thm:smoothbn2}. From the proof of Theorem~\ref{thm:smoothbn2}, we can represent $V_{-4}$ with the Kirby diagram in Figure~\ref{f:KirbyThmB13}, where we can see the $B_n \hookrightarrow V_{-4}$. If we rationally blow up this $B_n$, then we obtain the Kirby diagram in Figure~\ref{f:KirbyThmB12}, which by a sequence of Kirby moves gets us back to the diagram in Figure~\ref{f:KirbyThmB5}, which is precisely $V_{-4} \# (n-1)\overline{\C P^2}$.

\section{``Simple" embeddings}
\label{sec:simple}

The embeddings of the $B_n$s in Theorems~\ref{thm:smoothbn1} and \ref{thm:smoothbn2} are inherently different from the embeddings of $B_n \hookrightarrow E(m)_n$, as discussed in the beginning of section~\ref{sec:smoothemb}. As seen from Corollaries~\ref{cor:smoothbn1} and \ref{cor:smoothbn2}, the embeddings of $B_n \hookrightarrow V_{-n-1},V_{-4}$ are such that if one rationally blows up those $B_n$s and then performs the regular blow-down $(n-1)$ times, then one gets back the manifolds $V_{-n-1},V_{-4}$ respectively. One could also do these two steps in reverse: if one starts with $V_{-n-1},V_{-4}$, blows them up $(n-1)$ times and then rationally blows down the obtained $C_n$ configuration, then one again obtains the manifolds $V_{-n-1},V_{-4}$ respectively. This is summarized in the following diagrams for the embeddings of $B_n \hookrightarrow V_{-n-1}$ for $n \geq 2$ and for $B_n \hookrightarrow V_{-4}$ for odd $n \geq 3$ respectively:

$$
\begin{CD}
  V_{-n-1} @>\text{RBU the } B_n>> V_{-n-1}\# (n-1)\overline{\C P^2} \\
 @VV\text{BU }(n-1) \text{ times} V @VV\text{BD }(n-1) \text{ times}V  @.\\
  V_{-n-1}\# (n-1)\overline{\C P^2} @>\text{RBD the } C_n>> V_{-n-1} 
\end{CD}
$$

\vspace{.2in}

$$
\begin{CD}
  V_{-4} @>\text{RBU the } B_n>> V_{-4}\# (n-1)\overline{\C P^2} \\
 @VV\text{BU }(n-1) \text{ times} V @VV\text{BD }(n-1) \text{ times}V  @.\\
  V_{-4}\# (n-1)\overline{\C P^2} @>\text{RBD the } C_n>> V_{-4} 
\end{CD}
$$

\vspace{.2in}

This is not the case with the embeddings of $B_n \hookrightarrow E(m)_n$, since the rational blow-ups of those $B_n$'s result in $E(m) \# (n-1)\overline{\C P^2}$ (and not $E(m)_n \# (n-1)\overline{\C P^2}$) and so blowing down $(n-1)$ times yields the manifold $E(m)$ and not $E(m)_n$, the manifold we started with. 

As a result, one can call an embedding of $B_n \hookrightarrow X$ ``simple" if rationally blowing up and then blowing down $(n-1)$ times yields back the same 4-manifold $X$ (the top and right arrows of the diagram below). Equivalently, an embedding $B_n \hookrightarrow X$ is ``simple" if blowing up $(n-1)$ times followed by rationally blowing down the $C_n$, yields back the same 4-manifold $X$ (the left and bottom arrows of the diagram below).

$$
\begin{CD}
  X @>\text{RBU the } B_n>> X\# (n-1)\overline{\C P^2} \\
 @VV\text{BU }(n-1) \text{ times} V @VV\text{BD }(n-1) \text{ times}V  @.\\
  X\# (n-1)\overline{\C P^2} @>\text{RBD the } C_n>> X 
\end{CD}
$$

\vspace{.2in}
\noindent It follows that the embeddings of of $B_n \hookrightarrow V_{-n-1}$ for $n \geq 2$ and for $B_n \hookrightarrow V_{-4}$ for odd $n \geq 3$ are ``simple", whereas the embedding $B_n \hookrightarrow E(m)_n$ is not ``simple". Therefore, one can ask the following question: Are there obstructions to embedding the $B_n$s in a ``non-simple" way? 

Nevertheless, Theorem~\ref{thm:smoothbn2} prevents one from finding an upper bound on $n$ for a smooth $4$-manifold $X$ to contain an embedded $B_n$. However, one can ask whether such a bound exists for ``non-simple" embeddings of $B_n \hookrightarrow X$. 

The Kirby diagrams in the proofs of Theorems~\ref{thm:smoothbn1} and \ref{thm:smoothbn2} strongly suggest that the key to determining whether an embedding of a rational homology ball $B_n$ is ``simple" lies in analyzing how the extra sphere with self-intersection $(-1)$ intersects with the spheres of the $C_n$ configuration after one rationally blows up the $B_n$. For example, if one starts with $B_n \hookrightarrow V_{-n-1}$ for $n \geq 2$, and rationally blows it up, one obtains the Kirby diagram seen in Figure~\ref{f:KirbyThmA4}. In this case, the extra sphere with self-intersection $(-1)$ intersects with the last sphere of self-intersection $(-2)$ ($S_{n-1}$ in Figure~\ref{f:KirbyCnCircles1}) in the $C_n$ configuration. Likewise, if one starts with $B_n \hookrightarrow V_{-4}$ for odd $n \geq 3$, and rationally blows it up, one obtains the Kirby diagram seen in Figure~\ref{f:KirbyThmB5}. In this case, the extra sphere with self-intersection $(-1)$ intersects with the first sphere of self-intersection $(-2)$ ($S_{2}$ in Figure~\ref{f:KirbyCnCircles1}) in the $C_n$ configuration. In the ``non-simple" embedding case of $B_n \hookrightarrow E(m)_n$, if one rationally blows up those rational homology balls, then the extra sphere of self-intersection $(-1)$ intersects with the first and last spheres of the $C_n$ configuration ($S_{1}$ and $S_{n-1}$, respectively, in Figure~\ref{f:KirbyCnCircles1}), as seen in Figure~\ref{f:fishfiber}. 

In those instances where an exotic smooth manifold $X_{(n)}$ is obtained after rationally blowing down $X$ (with $B_n$), then the embedding $B_n \hookrightarrow X_{(n)}$ must be ``non-simple". Moreover, the intersection patterns of certain spheres of self-intersection $(-1)$ with the $C_n$ configuration back up in $X$, may be directly related to the obtained exotic smooth structure of $X_{(n)}$. Consequently, understanding the precise way of how the rational homology balls embed in $4$-manifolds, may give us a better understanding of exotic smooth structures of $4$-manifolds which were obtained as a result of a rational blow-down.


\end{document}